\documentclass[12pt]{article}

\usepackage{latexsym}
\usepackage{calc}
\usepackage{datetime}
\usepackage{amssymb, amsmath, amsfonts, listings, mathrsfs} 
\usepackage{xcolor} %from Anton
\usepackage[normalem]{ulem} %from Anton

\usepackage{graphicx}
\usepackage[labelfont=bf]{caption}

\newcounter{hours}\newcounter{minutes}

\newcommand{\printtime}{
\setcounter{hours}{\time/60}
\setcounter{minutes}{\time - \value{hours}*60}
\thehours h \theminutes m}

\newcommand{\stkout}[1]{\ifmmode\text{\sout{\ensuremath{#1}}}\else\sout{#1}\fi} %from Anton: \CDEL{Old Text $a = b$} \CADD{New Text $c = d$}.
\def\nr{\par \noindent}

\def\Diag{{\rm Diag \,}}
\def\Def{\stackrel{\mathrm{def}}{=}}

\def\inter{{\rm int \,}}

\def\wh{\widehat}
\def\ua{\underline{\alpha}}
\def\vf{\varphi}
\def\beq{\begin{equation}}
\def\eeq{\end{equation}}

\def\R{\mathbb{R}}

\def\BI{\begin{itemize}}
\def\EI{\end{itemize}}

\newcommand{\SetEQ}{\setcounter{equation}{0}}
\newcommand{\refLE}[1]{\ensuremath{\stackrel{(\ref{#1})}{\leq}}}
\newcommand{\refEQ}[1]{\ensuremath{\stackrel{(\ref{#1})}{=}}}
\newcommand{\refGE}[1]{\ensuremath{\stackrel{(\ref{#1})}{\geq}}}
\newcommand{\refGT}[1]{\ensuremath{\stackrel{(\ref{#1})}{>}}}

\newtheorem{theorem}{Theorem}
\newtheorem{lemma}{Lemma}
\newtheorem{corollary}{Corollary}

\newtheorem{assumption}{Assumption}
\newtheorem{definition}{Definition}

\newtheorem{example}{Example}
\newtheorem{remark}{Remark}
\newcommand{\proof}{\bf Proof: \rm \nr}
\newcommand{\qed}{\hfill $\Box$ \nr \medskip}
\newcommand{\half}{\mbox{${1 \over 2}$}}

\def\ba{\begin{array}}
\def\ea{\end{array}}
\def\beann{\begin{eqnarray*}}
\def\eeann{\end{eqnarray*}}
\def\bea{\begin{eqnarray}}
\def\eea{\end{eqnarray}}

\def\BT{\begin{theorem}}
\def\ET{\end{theorem}}
\def\BL{\begin{lemma}}
\def\EL{\end{lemma}}
\def\BC{\begin{corollary}}
\def\EC{\end{corollary}}
\def\BE{\begin{example}}
\def\EE{\end{example}}
\def\BD{\begin{definition}}
\def\ED{\end{definition}}
\def\BR{\begin{remark}}
\def\ER{\end{remark}}
\def\BAS{\begin{assumption}}
\def\EAS{\end{assumption}}
\def\BI{\begin{itemize}}
\def\EI{\end{itemize}}

\def\BMP{\begin{minipage}{9.5cm}}
\def\EMP{\end{minipage}}
\def\MPT{\begin{minipage}{11.5cm}}
\def\EPT{\end{minipage}}

\def\la{\langle}
\def\ra{\rangle}

\def\QF{\hspace{5ex} \Box}
\def\QR{\hfill \Box}

\title{Local and Global Convergence of Greedy Parabolic Target-Following Methods for Linear Programming}

\author{Yurii Nesterov\thanks{
CCOR at Corvinus Institute for Advanced Studies in Corvinus University of Budapest, and SDS in Chinese University of Hong Kong (Shenzhen). Professor emeritus at CORE (UCLouvain, Belgium). \newline Email: Yurii.Nesterov@uclouvain.be}
}

\date{December 17, 2024\\
\vspace{2ex} Printout: \printtime, \today
}

\begin{document}
\maketitle

\abstract{
In the first part of this paper, we prove that, under some natural non-degeneracy assumptions, the Greedy  Parabolic Target-Following Method \cite{PTS}, based on {\em universal tangent direction} \cite{AFF} has a favorable local behavior. In view of its global complexity bound of the order 
$O(\sqrt{n} \ln {1 \over \epsilon})$, 
this fact proves that the functional proximity measure, used for controlling the closeness to Greedy Central Path, is large enough for ensuring a local super-linear rate of convergence, provided that the proximity to the path is gradually reduced. 

This requirement is eliminated in our second algorithm, which is based on a new auto-correcting predictor direction. This method, besides the best-known polynomial-time complexity bound, ensures an automatic switching onto the local quadratic convergence in a small neighborhood of the solution. 

We present also the third algorithm, which approximates the path by quadratic curves. On the top of the best known global complexity bound, this method benefits from an unusual local  cubic rate of convergence. It is important that this amelioration needs no serious increase in the cost of one iteration.

Finally, we compare the advantages of these local accelerations with possibilities of finite termination. As we will see, the conditions allowing the optimal basis detection sometimes are even weaker than those required for the local superlinear convergence. Hence, it is important to endow the practical optimization schemes with both abilities. 

To the best of our knowledge, the proposed methods have a very interesting combination of favorable properties, which can be hardly found in the most of existing Interior-Point schemes. As all other parabolic target-following schemes, the new methods can start from an arbitrary strictly feasible primal-dual pair and go directly towards the optimal solution of the problem in a single phase. The preliminary computational experiments confirm the advantage of the second-order prediction. }

\vspace{1ex}

\section{Introduction}\label{sc-Intro}
\SetEQ

{\bf Motivation.} In the mid-eighties, starting from the seminal papers by Karmarkar \cite{Kar}, Renegar \cite{Jim}, and Gonzaga \cite{Clovis}, the Interior-Point Methods (IPM) for Linear Programming became the most active research direction in Optimization. The new methods, supported by very attractive worst-case polynomial-time complexity bounds, presented a serious competition for the traditional Simplex Method. Today, the most advanced versions of IPM are primal-dual predictor-corrector schemes, which follow the primal-dual central path in a large neighborhood, defined by some proximity measure (e.g. \cite{MTY}).

However, despite to the excellent complexity bounds, in the last years these methods are not very popular among practitioners. The reason is that the new problems of Machine Learning and Artificial Intelligence usually have a big dimension and a very special structure, which looks more suitable for the cheap first-order methods. However, the first-order methods are slow and suffer from the absence of polynomial-time complexity bounds. Hence, there is always a chance for adapting IPM to the new reality and getting even more efficient optimization schemes.

This paper presents one of the first steps in this direction. The main drawback of the classical theory of IPMs consists in the necessity of performing several stages of the minimization process (for explanations, see for example, Section 5.3.6 in \cite{Lect}). For the primal-dual pair of Linear Optimization Problems
\beq\label{prob-PDLO}
\ba{rcl}
\min\limits_{x \in \R^n_+} \{ \la c, x \ra: \; Ax=b \} & = & \max\limits_{y \in \R^m, s \in \R^n_+} \{ \la b, y \ra: \; s + A^T y = c \},
\ea
\eeq
the standard methods follows the {\em central path} $u_{\mu} = (x_{\mu},s_{\mu},y_{\mu})$ defined by the following system of equations:
$$
\ba{rcl}
x_{\mu} s_{\mu} & = & \mu e, \quad A x_{\mu} = b, \quad s_{\mu} + A^T y_{\mu} = c, \quad \mu > 0,
\ea
$$
where $e \in \R^n$ is the vector of all ones. Even if a feasible starting point
$$
\ba{rcl}
u_0 = (x_0, s_0 , y_0 ) & \in &{\cal F}_0 \Def \{ (x,s,y): \; Ax = b, \; s +A^T y = c, \; x, s \in \R^n_{++}\},
\ea
$$
is known,
we still need an initial stage for finding an approximation to the point $u_{\mu_0}$ for some $\mu_0>0$.

This stage can be eliminated in the framework of {\em weighted barriers} \cite{RTV}, where the {\em weighted central path} is defined by control variable $\bar v \in \R^n_{++}$ as follows:
\beq\label{def-WCP}
\ba{rcl}
x_{\bar v} s_{\bar v} & = & \bar v, \quad A x_{\bar v} = b, \quad s_{\bar v} + A^T y_{\bar v} = c, \quad \bar v > 0.
\ea
\eeq
However, it appears that the worst-case complexity bound in this approach depends on the condition number of the weights:
$$
\ba{rcl}
\kappa(\bar v) & = & \max\limits_{1 \leq i,j \leq n} {\bar v^{(i)} \over \bar v^{(j)}},
\ea
$$
and this destroys the polynomial-time complexity bounds of the schemes.

The latter difficulty was eliminated in \cite{PTS}, where the nonlinear equalities in (\ref{def-WCP}) were replaced by {\em convex inequalities}:
\beq\label{def-VCP}
\ba{rcl}
x s & \geq & v^2, \quad A x = b, \quad s + A^T y = c, \quad \la c, x \ra - \la b, y \ra \leq v^0,
\ea
\eeq
with $w=(v^0,v) \in \R^{n+1}$ being a vector of control parameters. The main advantage of this approach is related to the fact that the natural barrier function for the feasible set(\ref{def-VCP}),
$$
\ba{rcl}
F(u,w) & = & - \sum\limits_{i=1}^n \ln (x^{(i)} s^{(i)} - (v^{(i)})^2) - \ln (v^0 - \la c, x \ra + \la b, y \ra)
\ea
$$
admits a close-form solution for the problem
$$
\ba{rcl}
\min\limits_{u \in {\cal F}_0} F(u,w) & = & - (n+1)\ln{v^0 - \| v \|^2 \over n+1}.
\ea
$$
Thus, it is possible to measure the closeness of any point $u \in {\cal F}_0$ to the analytic center of the set (\ref{def-VCP}) by a simple {\em functional proximity measure} (FPM). The components of the control variable $w$ in this approach must satisfy the inequality $v^0 > \| v \|^2$, which explains the name {\em Parabolic Target Space}.

This idea was elaborated in \cite{PTS} in the framework of {\em self-concordant functions} (see, for example, Chapter 5 in \cite{Lect}). However, the corresponding machinery of Linear Algebra was quite heavy: instead of inverting at each iteration one $m \times m$-matrix, as it was done in the standard IPMs for Linear Optimization, it is necessary to invert $(2m)\times(2m)$-matrices.

This was the reason for revisiting this approach in the recent papers \cite{LCP,AFF}. In the second paper, we proposed a new {\em Universal Tangent Direction}, which is computationally cheap and which ensures the best-known worst-case complexity bound of $O(\sqrt{n} \ln {n \over \epsilon})$ for the number of Newton steps required for computing an $\epsilon$-solution of the problem (\ref{prob-PDLO}). The corresponding method can start from any point $u_0 \in {\cal F}_0$ and travel towards the optimal solution in a single stage.

In this paper, we start from further investigation of the properties of method \cite{AFF}. In particular, we prove for it a local linear convergence to the non-degenerate solution of the problem (\ref{prob-PDLO}) with coefficient depending only on the level of functional proximity level.

If this level vanishes, then we can get a super-linear convergence rate. However, a slight modification of the search direction gives us already a scheme with local quadratic convergence. Moreover, it is possible to replace the line-search strategy at the predictor step by a {\em parabolic search}. In this case, we get a local {\em cubic} rate of convergence. On the top of these results, we provide all our methods with a finite-termination criterion, which is based on the new indicator functions.

The classical results on local convergence and finite termination of IPMs for Linear Optimization are mainly based on Euclidean proximity measure \cite{Mer, Y1, YA, YY, ZTD, ZT, YZ}. Hence, our developments seem to be new. We confirm our theoretical results by  encouraging computational experiments, which confirm a superiority of the second-order prediction.

{\bf Contents.} The paper is organized as follows. In Section \ref{sc-Pred}, we introduce the framework of Parabolic Target Space \cite{PTS} and present a predictor-corrector method for the Greedy Strategy, based on the Universal Tangent Direction (UTD) \cite{AFF}. Our method~(\ref{met-AffinePT}) can be seen as a variant of Algorithm 4.1 in \cite{AFF}.

In Section \ref{sc-Loc}, under some natural non-degeneracy assumptions, we prove the local bounds for the size of some directions used in the method (\ref{met-AffinePT}). In Section \ref{sc-Par}, we derive a close-form expression for the growth of  FPM along UTD. It allows us to estimate the asymptotic local rate of convergence of the scheme, which appears to be linear with the 
coefficient $\half$.

In the next Section \ref{sc-Auto}, we define a new {\em auto-correcting} direction for the predictor step, which ensures the local quadratic convergence of the scheme. It also admits a worst-case global complexity bound of the order $O(\sqrt{n} \ln {n \over \epsilon})$.

Further, in Section \ref{sc-Second}, we define the second-order prediction strategy and prove for it the best global worst-case complexity bounds and the local {\em cubic} rate of convergence. The computational complexity of this scheme is essentially the same as that of the both previous schemes. However, as we will show in Section \ref{sc-Num} its computational behavior is much better.

Finally, in Section \ref{sc-Fin} we propose three new and easily computable indicators for finite termination of all our methods. 

{\bf Notations.} In this papers, the vectors in $\R^n$ are always denoted by lower-case Latin letters. An upper-case variant corresponds to the diagonal matrix:
$$
\ba{rcl}
x & \in \R^n, \quad X \; \Def \; \Diag(x) \in \R^{n \times n}.
\ea
$$
The positive orthant in $\R^n$ is denoted by $\R^n_+$ and for its interior we use notation $\R^n_{++}$.

For two vectors $x$ and $y$ of the same dimension, we denote by $\la x, y \ra$ its scalar product:
$$
\ba{rcl}
\la x, y \ra & = & \sum\limits_{i=1}^n x^{(i)} y^{(i)}, \quad x, y \in \R^n.
\ea
$$
We use the same notation $\la \cdot, \cdot \ra$ for vectors from different spaces. Hence, its actual sense is defined by the context. All arithmetic operations and relations involving vectors are understood in the component-wise sense. 

For Euclidean norm, we use notation
$$
\ba{rcl}
\| x \| & = & \la x, x \ra^{1/2}, \quad x \in \R^n.
\ea
$$
Similarly, $\ell_p$-norms with $p \geq 1$ are defined as follows:
$$
\ba{rcl}
\| x \|_p & = & \Big[ \sum\limits_{i=1}^n | x^{(i)}|^p \Big]^{1/p}, \quad x \in \R^n,
\ea
$$
with $\| x \|_{\infty} \Def \max\limits_{1 \leq i \leq n} | x^{(i)}|$. Note that for all $x,y \in \R^n$ we have
\beq\label{eq-Prod}
\ba{rcl}
\| x y \| & = & \Big[ \sum\limits_{i=1}^n \left( x^{(i)} y^{(i)}\right)^2 \Big]^{1/2} \; \leq \; \| x \|_4 \cdot \| y \|_4 \; \leq \; \| x \| \cdot \| y \|.
\ea
\eeq

For a matrix $C \in \R^{k \times p}$, we denote $\| C \|_{\infty} = \max\limits_{\forall (i,j)} | C^{(i,j)}|$. Then,
$$
\ba{rcl}
\| C x \|_{\infty} & \leq & \| C \|_{\infty} \| x \|_1, \quad \forall x \in \R^n.
\ea
$$

\section{Predictor-corrector scheme}\label{sc-Pred}
\SetEQ

Consider the standard primal-dual pair of Linear Programming problems:
\beq\label{prob-LP}
\ba{rcl}
\min\limits_{\ba{c} A x = b, \\ x \geq 0 \ea} \la c, x \ra & = & \max\limits_{\ba{c} s + A^T y = b, \\ s \geq 0 \ea} \la b, y \ra
\ea
\eeq
We assume existence of a strictly-feasible primal-dual solution $\hat u = (\hat x, \hat s, \hat y)$:
\beq\label{eq-SFeas}
\ba{rcl}
A \hat x & = & b, \quad \hat s + A^T \hat y \; = \; c, \quad \hat x, \hat s > 0.
\ea
\eeq
In what follows, we denote by ${\cal F}_0 = \Big\{ u = (x,s,y): Ax = b, s +A^T y = c, \; x,s \in \R^n_{++} \Big\}$ the relative interior of the feasible set of the primal-dual problem (\ref{prob-LP}).
For any $u \in {\cal F}_0$, we have the following useful relation:
\beq\label{eq-cost}
\ba{rcl}
\la c, x \ra - \la b, y \ra & = & \la s, x \ra.
\ea
\eeq

We solve the problem (\ref{prob-LP}) by the Parabolic Target Following approach \cite{PTS}, where the control variable $w = (v^{(0)},v) \in \R_+ \times \R^n$ is updated inside the {\em parabolic target set}
$$
\ba{rcl}
{\cal F}_p & = & \Big\{ w = (v^{(0)},v)  \in \R_+ \times \R^n: \; v^{(0)} > \| v \|^2 \Big\}.
\ea
$$
Sometimes we use notation $w^0 \Def v^{(0)}$ and $w^+ \Def v$.

For the barrier interpretation, let us introduce the full vector of variables $z=(u,w)$ belonging to the feasible set
$$
\ba{c}
{\cal F} = \Big\{ (u \in {\cal F}_0,w\in {\cal F}_p): \quad x^{(i)} s^{(i)} \geq (v^{(i)})^2, \; i = 1, \dots, n, \quad v^{(0)} \geq \la c, x \ra - \la b, y \ra \Big\}.
\ea
$$
This set admits a standard self-concordant barrier
$$
\ba{rcl}
F(z) & = & - \sum\limits_{i=1}^n \ln \left( x^{(i)} s^{(i)} - (v^{(i)})^2 \right) - \ln \left(v^{(0)} - \la c,x \ra + \la b, y \ra \right), \quad z \in \inter {\cal F}
\ea
$$
with parameter $\nu = 2n+1$. It can be shown \cite{PTS} that
$$
\ba{rcl}
\min\limits_{u: (u,w) \in {\cal F}} F(u,w) & = & \vf(w) \; = \; - (n+1) \ln \, \rho(w), \quad \rho(w) = \frac{v^{(0)} - \| v \|^2}{n+1}.
\ea
$$
Moreover, the optimal point $u(w) = (x(w),s(w),y(w))$ of this problem satisfies the following relations:
\beq\label{def-UW}
\ba{rcl}
x^{(i)}(w) s^{(i)}(w) & = &  (v^{(i)})^2 + \rho(w), \quad i=1, \dots, n,
\ea
\eeq
From these equations, we get
\beq\label{eq-SXW}
\ba{rcl}
\la s(w), x(w) \ra & = & \| v \|^2 + n \rho(w) \; = \; \frac{n v^{(0)} + \| v \|^2}{n+1}.
\ea
\eeq
Consequently,
\beq\label{eq-W0}
\ba{rcl}
v_0 - \la s(w), x(w) \ra & = & \rho(w).
\ea
\eeq

Note that the above relations justify the following {\em Functional Proximity Measure}:
\beq\label{def-FProx}
\ba{rcl}
\Psi(z) \; = \; F(z) - \vf(w) & = &  - \sum\limits_{i=1}^n \ln \left( x^{(i)} s^{(i)} - (v^{(i)})^2 \right) - \ln (v_0 - \la c,x \ra + \la b, y \ra)\\
& &  + (n+1) \ln \frac{v^{(0)} - \| v \|^2}{n+1} \; \geq \; 0,
\ea
\eeq
which vanishes only at points $z=(u(w),w)$ with $w \in {\cal F}_p$.

In our methods, we trace approximately the sequence $u(w_k)$ defined by the control variable $w_k \in {\cal F}_p$. The convergence $w_k \to 0$ is ensured by the simplest {\em Greedy Strategy}:
\beq\label{eq-Greedy}
\ba{rcl}
w_{k+1} & = & (1-\alpha_k) w_k, \quad \alpha_k \in (0,1), \; k \geq 0.
\ea
\eeq

Let us present an algorithmic description of our first method.
For its initialization, we need a strictly feasible point $u = (x,s,y) \in \inter {\cal F}_0$. By this point, we can define the control variable $w_*(u) = \left(v^{(0)}_*(u), v_*(u) \right)$ in the following way:
\beq\label{eq-Start}
\ba{rcl}
v_*^{(0)} & = & \la s, x \ra + \sigma(u), \quad v_*^{(i)}(w) \; = \; \sqrt{ x^{(i)} s^{(i)} - \sigma(u)}, \quad i = 1, \dots, n,
\ea
\eeq
where $\sigma(u) = \min\limits_{1 \leq i \leq n} x^{(i)} s^{(i)}$. It is easy to see that $u \refEQ{def-UW} u(w_*(u))$.

For an arbitrary pair $z = (u,w) \in {\cal F}$, in order to check closeness of $u$ to $u(w)$, we need to define the vector of residuals $r(z) \in \R^{n+1}$ as follows:
$$
\ba{rcl}
r^{(0)}(z) & = & v^{(0)} - \la s, x \ra \; \geq \; 0, \quad r^{(i)}(z) \; = \; x^{(i)} s^{(i)} - (v^{(i)})^2 \; \geq \; 0, \; i = 1, \dots, n.
\ea
$$
Note that
\beq\label{eq-Sum}
\ba{rcl}
\la r(z), e \ra & = & (n+1) \rho(w),
\ea
\eeq
where $e \in \R^{n+1}$ is the vector of all ones. Its truncated version is denoted by $\check e \in \R^n$. Similarly, vector $\check r(z) \in \R^n$ contains components of vector $r(z)$ with indexes $1 \leq i \leq n$.

We estimate the distances between points $u$ and $u(w)$ by the following measures:
\beq\label{def-Rho}
\ba{c}
\chi_k(z)\; = \; \left[ \sum\limits_{i=0}^n \frac{(r^{(i)}(z) - \rho(w))^2 }{ [r^{(i)}(z)]^k \, [\rho(w)]^{2-k}} \right]^{1/2}, \; k=0, 1, 2, \quad \delta(z) \; = \; \frac{\chi_1^2(z)}{ \chi_2(z)}.
\ea
\eeq
For $\chi_2(z) = 0$, define $\delta(z) = 0$.
If $r^{(i)}(z) = \rho(w)$ for all $0 \leq i \leq n$, then these measures vanish and $u=u(w)$ (see (\ref{def-UW}), (\ref{eq-W0})). Note that all these values are easy to compute. 

For the point $u_k=(x_k,s_k,y_k) \in {\cal F}$ and a right-hand side $d \in \R^n$, we define the {\em Universal Tangent  Direction} $\Delta_k(d) = (\Delta^x_k,\Delta^s_k,\Delta^y_k)(d)$ (see \cite{AFF}) as a unique solution of the following linear system:
\beq\label{def-ASK}
\ba{rcl}
X_k \Delta^s_k + S_k \Delta^x_k & = & d, \quad A \Delta^x_k = 0, \quad \Delta^s_k + A^T \Delta^y_k = 0,
\ea
\eeq
For its computation, we need to form and invert the matrix $\Sigma_k = A X_k S_k^{-1} A^T \in \R^{m \times m}$, which is independent on $d$.
We use also the following univariate function:
\beq\label{def-Theta*}
\ba{rcl}
\omega_*(\tau) & = &  - \tau - \ln(1-\tau), \quad 0 \leq \tau < 1.
\ea
\eeq
\beq\label{met-AffinePT}
\ba{|l|}
\hline \\
\hspace{2ex} \mbox{\bf Tangential Parabolic Target Following Method (TPTFM)}\\
\\
\hline \\
\mbox{{\bf Initialization.} Choose  $r \in (0,1)$, $A_{\psi} = \omega_*(r)$, $u_0 \in {\cal F}_0$, and $w_0 \refEQ{eq-Start} w_*(u_0)$.}\\
\mbox{Define the maximal proximity level $\beta = {r \over 2+r} <{1 \over 3}$.}\\
\\
\mbox{\bf $k$th iteration ($k \geq 0$).}\\
\\
\ba{rl}
\mbox{{\bf a)}} & \mbox{Compute $r(z_k)$ and $\Sigma_k^{-1} = \left[A X_k S_k^{-1} A^T\right]^{-1}$.}\\
& \mbox{Choose the acceptance level $\beta_k \in [0,\beta)$.}\\
\\
\mbox{{\bf b)}} & \mbox{If $\delta(z_k) \leq \beta_k$, then do \hspace{2ex} \fbox{\sc Predictor Step}}\\
& \bullet \; \mbox{Set $d_k = \left(\frac{\| v_k \|^2 }{ n+1} - \rho(w_k)\right) \check e - 2 v_k^2$ and compute $\Delta_k = \Delta_k(d_k)$.}\\
& \bullet \; \mbox{Define function  $\psi_k(\alpha) = \Psi(u_k+\alpha \Delta_k, (1-\alpha)w_k)$.}\\
& \bullet \; \mbox{Find $\alpha_k$ as an approximate solution of equation $\psi_k(\alpha) = A_{\psi}$.}\\
& \bullet \; \mbox{Define $u_{k+1} =u_k + \alpha_k \Delta_k$ and $w_{k+1} = (1-\alpha_k)w_k$.}\\
\\
\mbox{\bf c)} & \mbox{Otherwise, do \hspace{2ex} \fbox{\sc Corrector Step}}\\
& \bullet \; \mbox{Define $d_k = \rho(w_k) \check e - \check r(z_k)$. Compute $\Delta_k = \Delta_k(d_k)$.}\\
& \bullet \; \mbox{Define function  $f_k(\alpha) = F(u_k+\alpha \Delta_k, w_k)$.}\\
& \bullet \; \mbox{Find $\alpha_k$ as an approximate minimum of $f_k(\alpha)$ in $\alpha \geq 0$.}\\
& \bullet \; \mbox{Define $u_{k+1} =u_k + \alpha_k \Delta_k$ and $w_{k+1} = w_k$.}\\
\\
\mbox{\bf d)} & \mbox{If $w_k^0 \leq \epsilon$ and $\delta(z_k) \leq \beta_k$, then \fbox{\sc Stop}}
\ea\\
\\
\hline
\ea
\eeq
This method differs from the Algorithm 4.1 in \cite{AFF} mainly by a possibility to adjust the acceptance level $\beta_k \leq \beta$ during the minimization process. Our choice of $\beta$ ensures $r = {2 \beta \over 1 - \beta}$.

\section{Local size of Universal Tangent Direction}\label{sc-Loc}
\SetEQ

In this section, we justify the properties of the Universal Tangent Direction (\ref{def-ASK}) under the following non-degeneracy assumptions.
\BAS\label{ass-ND}
\begin{itemize}
\item
In problem (\ref{prob-LP}), there exists a unique primal solution $x^*$ with $m$ positive components. We assume that these are the first $m$ components of the vector:
$$
\ba{rcl}
x^* & = & (x^*_B, x^*_N), \quad x^*_B \in \R^m_{++}, \quad x^*_N = 0 \in \R^{n-m}.
\ea
$$
\item
In the corresponding partition $A = (A_B,A_N)$, the matrix $A_B \in \R^{m \times m}$ is non-degenerate. 
\item
Hence, $y^* = A_B^{-T} c_B$ (thus, $s^*_B = 0$), $x^*_B = A_B^{-1}b$, and we assume that $s^*_N \Def c_N - A_N^T y^* > 0$.
\end{itemize}
\EAS

From this assumption, we immediately derive several useful facts. Denote 
$$
\ba{c}
x^*_{\min} \; = \; \min\limits_{1 \leq i \leq m} x^*_i, \quad s^*_{\min} \; = \; \min\limits_{m+1 \leq i \leq m} s^*_i, \quad \pi_* \; = \; x^*_{\min} \cdot s^*_{\min},\\
\kappa \; = \; \| A^{-1}_B A_N \|_{\infty} \; = \; \| A_N^{T} A_B^{-T} \|_{\infty}.
\ea
$$
\BL\label{lm-Use}
Let $(x,s,y)$ be a feasible solution for the primal-dual problem (\ref{prob-LP}). Then, we have
\beq\label{eq-Gap}
\ba{rcl}
\la s^*_N, x_N \ra + \la s_B, x^*_B \ra & = & \la s, x \ra \; = \;  \la c, x \ra - \la b, y \ra,
\ea
\eeq
\beq\label{eq-DXYS}
\ba{rcl}
\| x_B - x^*_B \|_{\infty} & \leq & \kappa \| x_N \|_1, \quad
\| s_N - s^*_N \|_{\infty} \; \leq \; \kappa \| s_B \|_1.
\ea
\eeq
\EL
\proof Indeed,
$$
\ba{rcl}
0 & = & \la s - s^*, x - x^* \ra \; = \; \la s, x \ra - \la s^*,x \ra - \la s, x^* \ra,
\ea
$$
and we get (\ref{eq-Gap}). Further, from the definition of optimal partition, we have
$$
\ba{rcl}
A_B x^*_B & = & b \; = \; A_B x_B + A_N x_N,
\ea
$$
and we obtain the first inequality in (\ref{eq-DXYS}). Further, since
\beq\label{eq-RepDY}
\ba{rcl}
s_B & = & c_B - A_B^T y \; = \; A_B^T(y^*-y),
\ea
\eeq
we get
$$
\ba{rcl}
s_N - s_N^* & = & A_N^T (y^*-y) \; \refEQ{eq-RepDY} \; A_N^T  A_B^{-T} s_B,
\ea
$$
which results in the second inequality in (\ref{eq-DXYS}).
\qed
\BC\label{cor-Use}
Under conditions of Lemma \ref{lm-Use}, we have
\beq\label{eq-NonB}
\ba{rcl}
s^*_{\min} \| x_N \|_1 + x^*_{\min} \| s_B \|_1 & \leq & \la s, x \ra,
\ea
\eeq
\beq\label{eq-XSLow}
\ba{rcl}
x^{(i)} & \geq & {1 \over s^*_{\min}} ( \pi_* - \kappa \la s, x \ra), \; i = 1, \dots, m\\
\\
s^{(i)} & \geq & {1 \over x^*_{\min}} ( \pi_* - \kappa \la s, x \ra), \; i = m+1, \dots, n.
\ea
\eeq
\EC
\proof
Inequality (\ref{eq-NonB}) follows directly from (\ref{eq-Gap}). The first inequality in (\ref{eq-XSLow}) can be obtained as follows:
$$
\ba{rcl}
x^{(i)} & \geq & x^*_{\min} - \| x - x^*\|_{\infty} \; \refGE{eq-DXYS} \; x^*_{\min} - \kappa \| x_N \|_1 \\
\\
& \refGE{eq-NonB} & x^*_{\min} - {\kappa \over s^*_{\min}} \la s, x \ra \; = \; {1 \over s^*_{\min}} ( \pi_* - \kappa \la s, x \ra).
\ea
$$
The second inequality can be justified in the same way.
\qed

Let us apply Lemma \ref{lm-Use} for estimating the size of Universal Tangent Direction (UTD), defined by some positive definite diagonal matrices $X$ and $S$ and the following system of linear equations:
\beq\label{eq-LSys}
\ba{rcl}
S \Delta^x + X \Delta^s & = & d, \quad A \Delta^x \; = \; 0, \quad \Delta^s + A^T \Delta^y \; = \; 0,
\ea
\eeq
with some $d \in \R^n$. Denote
\beq\label{def-Norms}
\ba{rcl}
\delta_x & = & \| X_B^{-1} d_B \|, \quad \rho_x = \| X^{-1}_B S_B \| \; = \; \max\limits_{1 \leq i \leq m} {s^{(i)} \over x^{(i)}}, \\
\\
\delta_s & = & \| S_N^{-1}d_N \|, \quad \rho_s \; = \; \| S_N^{-1} X_N \| \; = \; \max\limits_{m+1 \leq i \leq n} {x^{(i)} \over s^{(i)}}.
\ea
\eeq

\BT
Let the feasible primal-dual point $(x,s,y)$ is close enough to the optimal solution:
\beq\label{eq-CCond}
\ba{rcl}
\rho_x \rho_s & < & \kappa^{-2}.
\ea
\eeq
Then the size of the UTD (\ref{eq-LSys}) is bounded as follows:
\beq\label{eq-DBound}
\ba{rcl}
\| \Delta^x_B \| & \leq & {\kappa \over 1 - \kappa^2 \rho_x \rho_s} [\delta_s + \kappa  \rho_s \delta_x], \quad \| \Delta^x_N \| \; \leq \; {1 \over 1 - \kappa^2 \rho_x \rho_s} [\delta_s + \kappa  \rho_s \delta_x],\\
\\
\| \Delta^s_B \| & \leq & {1 \over 1 - \kappa^2 \rho_x \rho_s} [\delta_x + \kappa  \rho_x \delta_s], \quad
\| \Delta^s_N \| \; \leq \; {\kappa \over 1 - \kappa^2 \rho_x \rho_s}
[\delta_x + \kappa  \rho_x \delta_s].
\ea
\eeq
\ET
\proof
Let us represent the solution of the system (\ref{eq-LSys}) in terms of the optimal partition. Note that
$$
\ba{rcl}
\Delta_B^s & = & X^{-1}_B (d_B - S_B \Delta_B^x), \quad \Delta^x_B \; = \; - A_B^{-1} A_N \Delta^x_N, \quad \Delta^x_N \; = \; S_N^{-1}(d_N - X_N \Delta^s_N).
\ea
$$
Hence,
$$
\ba{rcl}
\Delta_B^s & = & X^{-1}_B \left(d_B + S_B A_B^{-1} A_N S_N^{-1}(d_N - X_N \Delta^s_N) \right).
\ea
$$
At the same time, $\Delta^s_B = - A^T_B \Delta^y$ and $\Delta^s_N = - A^T_N \Delta_y$. Hence,
$\Delta^s_N = A^T_N A_B^{-T} \Delta^s_B$,
and we conclude that
$$
\ba{rcl}
\Delta_B^s & = & X^{-1}_B \left(d_B + S_B A_B^{-1} A_N S_N^{-1}(d_N - X_N A^T_N A_B^{-T} \Delta^s_B) \right)\\
\\
& = & X^{-1}_B d_B + X^{-1}_B S_B A_B^{-1} A_N S_N^{-1}d_N - X^{-1}_B S_B A_B^{-1} A_N S_N^{-1} X_N A^T_N A_B^{-T} \Delta^s_B.
\ea
$$
Then for $d$ small enough, by the representation above, we get
$$
\ba{rcl}
\| \Delta^s_B \| & \leq & {1 \over 1 - \kappa^2 \rho_x \rho_s} [\delta_x + \kappa  \rho_x \delta_s].
\ea
$$
At the same time,
$$
\ba{rcl}
\Delta^s_N & = & A^T_N A_B^{-T} X^{-1}_B \left(d_B + S_B A_B^{-1} A_N S_N^{-1}(d_N - X_N\Delta^s_N) \right) \; = \;  A^T_N A_B^{-T} X^{-1}_B d_B \\
\\
& & + A^T_N A_B^{-T} X^{-1}_B S_B A_B^{-1} A_N S_N^{-1}d_N -  A^T_N A_B^{-T} X^{-1}_B S_B A_B^{-1} A_N S_N^{-1} X_N\Delta^s_N,
\ea
$$
and we conclude that
$$
\ba{rcl}
\| \Delta^s_N \| & \leq & {\kappa \over 1 - \kappa^2 \rho_x \rho_s}
[\delta_x + \kappa  \rho_x \delta_s].
\ea
$$

The remaining inequalities can be obtained by the following representations:
$$
\ba{rcl}
\Delta_B^x & = & - A_B^{-1} A_N S_N^{-1} ( d_N - X_N A_N^T A_B^{-T} X_B^{-1}( d_B - S_B \Delta_B^x)),\\
\\
\Delta_N^x & = & S_N^{-1} ( d_N - X_N A_N^T A_B^{-T} X_B^{-1}( d_B + S_B A_B^{-1} A_N \Delta_N^x)). \QF
\ea
$$

We need some sufficient conditions for inequality (\ref{eq-CCond}).
\BL\label{lm-Suff}
Let the feasible primal-dual point $(x,s,y)$ be close enough to the optimal solution:
\beq\label{eq-CLC}
\ba{rcl}
\la s, x \ra & < & {\pi_* \over \kappa}.
\ea
\eeq
Then we have the following bounds:
\beq\label{eq-RBound}
\ba{rcl}
\rho_x & \leq & {s^*_{\min} \la s, x \ra \over x^*_{\min} (\pi_* - \kappa \la s, x \ra)}, \quad \rho_s \; \leq \; {x^*_{\min} \la s, x \ra \over s^*_{\min} (\pi_* - \kappa \la s, x \ra)}, 
\ea
\eeq
\beq\label{eq-DBound1}
\ba{rcl}
\delta_x &\leq &  {s^*_{\min} \| d_B \| \over \pi_* - \kappa \la s, x \ra}, \quad \delta_s \; \leq \; {x^*_{\min} \| d_N \| \over \pi_* - \kappa \la s, x \ra}.
\ea
\eeq
\EL
\proof
Indeed, $\rho_x = \| X_B^{-1} S_B \| \refLE{eq-XSLow} {s^*_{\min} \over \pi_* - \kappa \la s, x \ra} \| S_B \| \refLE{eq-NonB} {s^*_{\min} \over \pi_* - \kappa \la s, x \ra} \cdot {\la s, x \ra \over x^*_{\min}}$. The second inequality in~(\ref{eq-RBound}) can be proved in a similar way. 
The remaining inequalities in (\ref{eq-DBound1}) also follow from (\ref{eq-XSLow}).
\qed

Let us specify the upper bounds (\ref{eq-DBound}) in the following  neighbourhood of the solution:
\beq\label{eq-Neib}
\ba{rcl}
\la s, x \ra & \leq & {\pi_* \over 4 \kappa}.
\ea
\eeq 
\BL\label{lm-DBound2}
Let the feasible primal-dual pair $(x,s,y)$ satisfy condition (\ref{eq-Neib}). Then
\beq\label{eq-DBound2}
\ba{rcl}
\| \Delta^x_B \| \cdot \|\Delta^s_B \| & \leq & {2\kappa \over \pi_*}\| d \|^2 , \quad
\| \Delta^x_N \| \cdot \|\Delta^s_N \| \; \leq \;  {2\kappa \over \pi_*}\| d \|^2.
\ea
\eeq
Moreover,
\beq\label{eq-TBound}
\ba{rcl}
\| \Delta^x \| & \leq &  \sqrt{{5 \over 2}(1+\kappa^2)} \, {\| d \|\over s^*_{\min}} , \quad \| \Delta^s \| \; \leq \;\sqrt{{5 \over 2}(1+\kappa^2)} \, {\| d \|\over x^*_{\min}} .
\ea
\eeq
\EL
\proof
Denote $\epsilon = \la s, x \ra$.
Then, in view of inequalities (\ref{eq-RBound}), we have
$\kappa^2 \rho_x \rho_s \leq { \kappa^2 \epsilon^2 \over (\pi_* - \kappa \epsilon)^2} \refLE{eq-Neib} {1 \over 9}$. At the same time,
$$
\ba{rcl}
\delta_s + \kappa \rho_s \delta_x & \stackrel{(\ref{eq-RBound}),(\ref{eq-DBound1})}{\leq} & {x^*_{\min} \| d_N \| \over \pi_* - \kappa \epsilon} + \kappa {x^*_{\min} \epsilon \over s^*_{\min} (\pi_* - \kappa \epsilon)} \cdot {s^*_{\min} \| d_B \| \over \pi_* - \kappa \epsilon}\\
\\
& = & {x^*_{\min}  \over \pi_* - \kappa \epsilon} \left[ \| d_N \| +{\kappa \epsilon \| d_B \| \over \pi_* - \kappa \epsilon} \right] \; \refLE{eq-Neib} \; {x^*_{\min}  \over \pi_* - \kappa \epsilon} \left[ \| d_N \| +{1 \over 3} \| d_B \| \right].
\ea
$$
Thus, $\| \Delta_B^x \| \refLE{eq-DBound} {9 \kappa x^*_{\min}  \over 8(\pi_* - \kappa \epsilon)} \left[ \| d_N \| +{1 \over 3} \| d_B \| \right]$. Similarly,
$$
\ba{rcl}
\delta_x + \kappa \rho_x \delta_s & \stackrel{(\ref{eq-RBound}),(\ref{eq-DBound1})}{\leq} & {s^*_{\min} \| d_B \| \over \pi_* - \kappa \epsilon} + \kappa {s^*_{\min} \epsilon \over x^*_{\min} (\pi_* - \kappa \epsilon)} \cdot {x^*_{\min} \| d_N \| \over \pi_* - \kappa \epsilon}\\
\\
& = & {s^*_{\min}  \over \pi_* - \kappa \epsilon} \left[ \| d_B \| +{\kappa \epsilon \| d_N \| \over \pi_* - \kappa \epsilon} \right] \; \refLE{eq-Neib} \; {s^*_{\min}  \over \pi_* - \kappa \epsilon} \left[ \| d_B \| +{1 \over 3} \| d_N \| \right].
\ea
$$
Thus, $\| \Delta_B^s \| \refLE{eq-DBound} {9 s^*_{\min}  \over 8(\pi_* - \kappa \epsilon)} \left[ \| d_B \| +{1 \over 3} \| d_N \| \right]$. Note that for two numbers $a,b \geq 0$, we have 
$$
\ba{rcl}
(a + {1 \over 3}b)({1 \over 3} a + b) & \leq & {8 \over 9}(a^2+b^2).
\ea
$$
Hence, we conclude that
$$
\ba{rcl}
\| \Delta^x_B \| \cdot \| \Delta^s_B \| & \leq & {9^2 \kappa \pi_*   \over 8^2(\pi_* - \kappa \epsilon)^2} \left[ \| d_B \| +{1 \over 3} \| d_N \| \right] \cdot \left[ {1 \over 3}\| d_B \| + \| d_N \| \right] \\
\\
& \leq & {9 \kappa \pi_*   \over 8(\pi_* - \kappa \epsilon)^2} \left[ \| d_B \|^2 +\| d_N \|^2 \right] \; \refLE{eq-Neib} \; {2\kappa \over \pi_*} \left[ \| d_B \|^2 +\| d_N \|^2 \right].
\ea
$$

Similarly, since
$$
\ba{rcl}
\delta_s + \kappa \rho_s \delta_x & \leq & {x^*_{\min}  \over \pi_* - \kappa \epsilon} \left[ \| d_N \| +{1 \over 3} \| d_B \| \right], \quad
\delta_x + \kappa \rho_x \delta_s \; \leq \; {s^*_{\min}  \over \pi_* - \kappa \epsilon} \left[ \| d_B \| +{1 \over 3} \| d_N \| \right],
\ea
$$
we have
$$
\ba{rcl}
\| \Delta^x_N \| \cdot \| \Delta^s_N \| & \leq & {9^2 \kappa \pi_*   \over 8^2(\pi_* - \kappa \epsilon)^2} \left[ \| d_N \| +{1 \over 3} \| d_B \| \right] \cdot \left[ {1 \over 3}\| d_N \| + \| d_B \| \right] \\
\\
& \leq & {9 \kappa \pi_*   \over 8(\pi_* - \kappa \epsilon)^2} \left[ \| d_B \|^2 +\| d_N \|^2 \right] \; \refLE{eq-Neib} \; {2\kappa \over \pi_*} \left[ \| d_B \|^2 +\| d_N \|^2 \right]. 
\ea
$$

Finally, in view of (\ref{eq-Neib}) and relation $(a+{1 \over 3}b)^2 \leq {10 \over 9}(a^2+b^2)$, we have
$$
\ba{rcl}
\| \Delta^x \|^2 & = & \| \Delta^x_B\|^2 + \| \Delta^x_N \|^2 \; \refLE{eq-DBound} \; \left({9 \over 8}\right)^2 (1+\kappa^2) (\delta_s + \kappa \rho_s \delta_x)^2 \\
\\
& \leq & \left({9 \over 8}\right)^2 (1+\kappa^2) \left( {4 x^*_{\min} \over 3 \pi_*}\Big[ \| d_B \| + {1 \over 3}\| d_N \|\Big] \right)^2\; \leq \; {5 \over 2} (1+\kappa^2)\left( {1 \over s^*_{\min}}\right)^2 \| d \|^2.
\ea
$$
The second inequality in (\ref{eq-TBound}) can be proved in a similar way.
\qed

The statement of Lemma \ref{lm-DBound2} leads to the following important consequence:
\beq\label{eq-SumBND}
\ba{rcl}
\| \Delta^x \Delta^s \|^2 & = & \| \Delta^x_B \Delta^s_B \|^2 + 
\| \Delta^x_N \Delta^s_N \|^2 \; \leq \; 
\| \Delta^x_B \Delta^s_B \|^2_1 + \| \Delta^x_N \Delta^s_N \|^2_1
\\ \\
& \leq & \| \Delta^x_B \|^2 \cdot \| \Delta_B^s \|^2 + \| \Delta^x_N \|^2 \cdot \| \Delta^s_N \|^2 \; \refLE{eq-DBound2} \; {8 \kappa^2\over \pi_*^2} \| d \|^4.
\ea
\eeq

\section{Local predictor abilities of TPTFM}\label{sc-Par}
\SetEQ

Let us estimate the performance of method (\ref{met-AffinePT}) at the predictor step. For this regime, we have
$$
\ba{rcl}
\delta(z_k) & \leq & \beta_k.
\ea
$$
In accordance to Lemma 5.3 in \cite{AFF} and inequalities (5.10), (5.11) there, this implies the following relations:
\beq\label{eq-Prop}
\ba{rcl}
(1-\beta_k)\rho(w_k) e & \leq & r(z_k) \; \leq \; {1 \over 1- \beta_k} \rho(w_k) e, 
\ea
\eeq
\beq\label{eq-PropXS}
\ba{rcl}
(1-\beta_k)x(w_k)s(w_k) & \leq & x s \; \leq \; {1 \over 1- \beta_k} x(w_k)s(w_k).
\ea
\eeq
Moreover, we have: \footnote{$^)$ In \cite{AFF}, the first inequality in (\ref{eq-Rho}) was obtained inside the proof of Lemma 5.3 in the form $\zeta_0(z) \leq {\beta \over \sqrt{1-\beta}}$. The second inequality in (\ref{eq-Rho}) was obtained in the proof of Lemma 5.4 as the relation (5.13).}$^)$
\beq\label{eq-Rho}
\ba{rcl}
\chi_1(z_k) & \leq & {\beta_k \over \sqrt{1-\beta_k}}, \quad \chi_0(z_k) \; \Def \; {1 \over \rho(w_k)} \| r(z_k) - \rho(w_k)  e \| \; \leq \; {\beta_k \over 1 - \beta_k}. 
\ea
\eeq

For the sake of notation, let us drop index $k$ for all objects related to the $k$th iteration.
By equality (5.16) in \cite{AFF}, for the predictor step $z(\alpha) = z + \alpha \Delta$, we have
\beq\label{eq-PGrow}
\ba{rcl}
\Psi(z(\alpha)) & = & - \sum\limits_{i=0}^n \ln \left( 1 + {1 \over \rho(w(\alpha))} A^{(i)}(\alpha) \right),
\ea
\eeq
where 
\beq\label{eq-OLow}
\ba{rcl}
\rho(w(\alpha)) & = & {1 \over n+1} \Big((1-\alpha)v_0 - (1-\alpha)^2 \| v \|^2\Big) \; \geq \; (1-\alpha)  \rho(w),
\ea
\eeq
and $A(\alpha) = r(z) - \rho(w) e + \alpha^2 g \in \R^{n+1}$ with $g^{(0)} = {1 \over n+1} \| v \|^2$, and
\beq\label{def-GG}
\ba{rcl}
g^{(i)} & = & (\Delta^x \Delta^s)^{(i)} - (v^{(i)})^2 + {1 \over n+1} \| v \|^2, \quad i = 1, \dots, n.
\ea
\eeq
Note that $\la A(\alpha), e \ra = 0$. If ${1 \over \rho(w(\alpha))} \| A(\alpha) \| < r$, then by the rules of the method, we have
$$
\ba{rcl}
\omega_*\left( r \right) & = & A_{\psi} \; = \; 
\Psi(z(\alpha)) \;  <  \; \omega_*\left( r \right),
\ea
$$
and this is impossible. Hence, since $\alpha < 1$, we conclude that
\beq\label{eq-Size}
\ba{c}
r \leq {1 \over \rho(w(\alpha))} \| A(\alpha) \| \leq {1 \over \rho(w(\alpha))} \left[ \| r(z) - \rho(w) e\| + \| \Delta^x \Delta^s \| + \| {1 \over n+1} \| v \|^2  e  - v^2_+ \| \right],
\ea
\eeq
where $v_+ = (0,v) \in \R^{n+1}$. Let us estimate separately the terms in the right-hand side. We have
$$
\ba{rcl}
{1 \over \rho(w(\alpha))} \| r(z) -\rho(w) e\| & = & {\rho(w) \over \rho(w(\alpha))} \chi_0(z) \; \stackrel{(\ref{eq-Rho}),(\ref{eq-OLow})}{\leq} \; {\beta_k \over (1 - \alpha)(1-\beta_k)}.
\ea
$$
Further, in view of inequality (\ref{eq-SumBND}), we have
$\| \Delta^x \Delta^s \| \leq 2 \sqrt{2} {\kappa \over \pi_*} \| d_a \|^2$, where $d_a$ is the right-hand side applied in Item b) of (\ref{met-AffinePT}). Then,
\beq\label{eq-DA}
\ba{c}
\| d_a \|^2 \; = \; \| (\rho(w) - {1 \over n+1} \| v \|^2) \check e + 2 v^2 \|^2 \; = \; \| 2 (v^2 + \rho(w) \check e) - {v^{(0)} \over n+1} \check e \|^2\\
\\
\; = \; 4 \sum\limits_{i=1}^n \left( (v^{(i)})^2 + \rho(w) \right)^2 - {4 v^{(0)}  \over n+1}\left(\| v \|^2 + n \rho(w) \right) + {n (v^{(0)})^2  \over (n+1)^2} \\
\\
\; \leq \; 4 \| v \|^4 + 8 \rho(w) \| v \|^2 + 4 n \rho^2(w) - {4 (v^{(0)})^2 \over n+1}\left(\| v \|^2 + n \rho(w) \right) + {n (v^{(0)})^2  \over (n+1)^2}\\
\\
\; = \; 4 \| v \|^4 + 8 \rho(w) \| v \|^2 + 4 n \rho^2(w) - 4 \rho(w) \left(\| v \|^2 + n \rho(w) \right) \\
\\
- {4 \| v \|^2 \over n+1}\left(\| v \|^2 + n \rho(w) \right) + {n (v^{(0)})^2  \over (n+1)^2}\\
\\
\; = \; {n (v^{(0)})^2 \over (n+1)^2} +  {4 \rho(w) \| v \|^2 \over n+1}  + {4n \| v \|^4 \over n+1} \; \leq \; {(v^{(0)})^2 \over n+1} +  {4 v^{(0)} \| v \|^2 \over (n+1)^2}  + 4\| v \|^4 .
\ea
\eeq
Finally,
\beq\label{eq-DV4}
\ba{rcl}
\| {1 \over n+1} \| v \|^2   e  - v^2_+ \|^2 & = & {\| v \|^4 \over n+1} - {2 \over n+1}\| v \|^4 +  \| v \|^4_4 \; \leq \; \| v \|^4.
\ea
\eeq

Thus, we have proved the following bound:
\beq\label{eq-BB}
\ba{rcl}
r & \leq & {\beta \over (1 - \alpha)(1-\beta)} + {n+1 \over (1-\alpha)(v^{(0)} - \| v \|^2)} \left[ \| v \|^2 +  {2^{3/2}\kappa \over \pi_*} \left( {(v^{(0)})^2 \over n+1} +  {4 v^{(0)} \| v \|^2 \over (n+1)^2}  + 4\| v \|^4 \right) \right].
\ea
\eeq

Let us look now at the behavior of method (\ref{met-AffinePT}) from the global perspective. Note that all control variables in this scheme have the following representation:
$$
\ba{rcl}

w_0 \; \Def \; (\bar v_0, \bar v), \quad w_k & = & \tau_k w_0, \quad \tau_0 = 1, \quad \tau_{k+1} = (1-\alpha_k) \tau_k, \quad k \geq 0.
\ea
$$
Denoting $f_k \Def v^{(0)}_k = \tau_k \bar v_0$, we can rewrite (\ref{eq-BB}) as follows:
\beq\label{eq-BB1}
\ba{rcl}
(1 - \alpha_k)r & \leq & {\beta_k \over 1-\beta_k} + {(n+1) \tau_k^2  \over \tau_k(\bar v_0 - \tau_k \| \bar v \|^2)} \left[ \| \bar v \|^2 + {2^{3/2}\kappa \over \pi_*} \left( {\bar v_0^2 \over n+1} +  {4 \bar v_0 \| \bar v \|^2 \tau_k \over (n+1)^2}  + 4 \tau_k^2 \| \bar v \|^4 \right)\right]\\
\\
& \leq & {\beta_k \over 1-\beta_k} + {(n+1) \tau_k  \over f_0 - f_{k}} \left[ \bar v_0 + {2^{3/2}\kappa \over \pi_*} \left( {\bar v_0^2 \over n+1} +  {4 \bar v^2_0 \tau_k  \over (n+1)^2}  + 4 \tau_k^2 \bar v_0^2 \right)\right]\\
\\
& = & {\beta_k \over 1-\beta_k} + {(n+1) f_k  \over f_0 - f_{k}} \left[ 1 + {2^{3/2}\kappa \over \pi_*} \left( {f_0 \over n+1} +  {4 f_k  \over (n+1)^2}  + 4 {f_k^2 \over f_0} \right)\right].
\ea
\eeq
Hence, since $f_k \leq f_0$, we have
$$
\ba{rcl}
f_{k+1} & = & (1-\alpha_k) f_k \; \leq \; {1 \over r}  \left\{ {\beta_k \over 1-\beta_k} + {(n+1) f_k  \over f_0 - f_{k}} \left[ 1 + {2^{3/2}\kappa \over \pi_*}  f_0 \left( {1 \over n+1} +  {4  \over (n+1)^2}  + 4 \right)\right]\right\} f_k\\
\\
& \leq & {1 \over r}  \left\{ {\beta_k \over 1-\beta_k} + {(n+1) f_k  \over f_0 - f_{k}} \left[ 1 + 11 \sqrt{2} {\kappa \over \pi_*} f_0 \right]\right\} f_k.
\ea
$$
Note that our bounds are valid only locally, when $\la s, x \ra \leq {\pi_* \over 4 \kappa}$. Thus, we have proved the following statement.
\BT\label{th-Quad}
Let $f_0 \leq {\pi_* \over 4 \kappa}$. Then, in the method (\ref{met-AffinePT}), we have
\beq\label{eq-QRate}
\ba{rcl}
f_{k+1} & \leq & {1 \over r}  \left\{ {\beta_k \over 1-\beta_k} + 5{(n+1) f_k  \over (f_0 - f_{k})}\right\} f_k
, \quad k \geq 0. 
\ea
\eeq
\ET

Our reasoning demonstrates a certain advantage of tracing the classical central path. In this case, $\bar v = 0$. Consequently,
$$
\ba{rcl}
r & \leq & {\beta_k \over (1-\alpha_k)(1-\beta_k)} + { f_k \over 1-\alpha_k} \cdot 2\sqrt{2} {\kappa \over \pi_*},
\ea
$$
and we have
$$
\ba{rcl}
f_{k+1} & =  & (1-\alpha_k) f_k \; \leq\; {1 \over r} \left[ {\beta_k \over 1-\beta_k} + 2\sqrt{2}{\kappa \over \pi_*}  f_k  \right]f_k
\; \leq \; {1 \over r} \left[ {\beta_k \over 1-\beta_k} + {f_k \over \sqrt{2} f_0} \right]f_k .
\ea
$$

From our analysis, we conclude that for the local quadratic convergence, we need to choose $\beta_k$ proportionally to $f_k$. If we keep $\beta_k$ constant, then the asymptotic convergence rate is linear, with the coefficient being an absolute constant. 
For example, for the choice $\beta_k = \beta$, we have ${\beta \over r(1-\beta)} = \half$, and the local rate is $\left(\half\right)^k$. 
Many other strategies for relating $\beta_k$ with $f_k$ are possible. However, in the remaining part of the paper, we will try to avoid these complications by improving the search direction at the predictor step.

\section{Auto-correcting predictor step}\label{sc-Auto}
\SetEQ

The main drawback of method (\ref{met-AffinePT}) is related to the fact that during its predictor step, the initial proximity measure can only increase.  
If the acceptance level $\beta_k$ is constant, this feature prevents the local quadratic convergence of the scheme (see (\ref{eq-QRate})). In this section, we analyze another version, where for the predictor Step b) we use a new right-hand side
$$
\ba{rcl}
\tilde d_k & = & {v_k^{(0)} \over n+1} \check e - x(w_k)s(w_k) - x_k s_k \; \refEQ{def-UW} \; \frac{\| v_k \|^2 }{ n+1} \check e - v_k^2 - x_k s_k.
\ea
$$
\beq\label{met-AutoPT}
\ba{|l|}
\hline \\
\hspace{2ex} \mbox{\bf Auto-Correcting Parabolic Target Following Method (ACPTFM)} \quad \\
\\
\hline \\
\mbox{{\bf Initialization.} Choose  $r \in (0,1)$, $A_{\psi} = \omega_*(r)$, $u_0 \in {\cal F}_0$, and $w_0 \refEQ{eq-Start} w_*(u_0)$.}\\
\mbox{Define the acceptance level $\beta = {r \over 2+r} < {1 \over 3}$.}\\
\\
\mbox{\bf $k$th iteration ($k \geq 0$).}\\
\\
\ba{rl}
\mbox{{\bf a)}} & \mbox{Compute $r(z_k)$ and $\Sigma_k^{-1} = \left[A X_k S_k^{-1} A^T\right]^{-1}$. }\\
\\
\mbox{{\bf b)}} & \mbox{If $\delta(z_k) \leq \beta$, then do \hspace{2ex} \fbox{\sc Predictor Step}}\\
& \bullet \; \mbox{Set $\tilde d_k = \frac{\| v_k \|^2 }{ n+1} \check e - v_k^2 - x_k s_k$ and compute $\tilde \Delta_k = \Delta_k(\tilde d_k)$.}\\
& \bullet \; \mbox{Define function  $\psi_k(\alpha) = \Psi(u_k+\alpha \tilde \Delta_k, (1-\alpha)w_k)$.}\\
& \bullet \; \mbox{Find $\alpha_k$ as an approximate solution of equation $\psi_k(\alpha) = A_{\psi}$.}\\
& \bullet \; \mbox{Define $u_{k+1} =u_k + \alpha_k \tilde \Delta_k$ and $w_{k+1} = (1-\alpha_k)w_k$.}\\
\\
\mbox{\bf c)} & \mbox{Otherwise, do \hspace{2ex} \fbox{\sc Corrector Step}}\\
& \bullet \; \mbox{Define $d_k = \rho(w_k) \check e - \check r(z_k)$. Compute $\Delta_k = \Delta_k(d_k)$.} \\
& \bullet \; \mbox{Define function  $f_k(\alpha) = F(u_k+\alpha \Delta_k, w_k)$.}\\
& \bullet \; \mbox{Find $\alpha_k$ as an approximate minimum of $f_k(\alpha)$ in $\alpha \geq 0$.}\\
& \bullet \; \mbox{Define $u_{k+1} =u_k + \alpha_k \Delta_k$ and $w_{k+1} = w_k$.}\\
\\
\mbox{\bf d)} & \mbox{If $w_k^{(0)} \leq \epsilon$ and $\delta(z_k) \leq \beta$, then \fbox{\sc Stop}}
\ea\\
\\
\hline
\ea
\eeq

At point $z=(u,w) \in {\cal F}$, let us define the right-hand sides $d_a$ as at Step b) and $d_c$ as at Step c) of method (\ref{met-AffinePT}). Then, the right-hand side $\tilde d$ of Step b) in method (\ref{met-AutoPT}), can be seen as a combination of these two vectors:
\beq\label{eq-DTSum}
\ba{rcl}
\tilde d & = & d_a + d_c \; = \; 
\frac{\| v \|^2 }{ n+1} \check e - v^2 - x s.
\ea
\eeq
As before, for points $\tilde z(\alpha) = z + \alpha (\tilde \Delta, - w)$ with $\alpha \geq 0$ and $\tilde \Delta = \Delta(\tilde d)$ (see (\ref{eq-LSys})), we can derive a closed-form expression for the values of functional proximity measure. 

Indeed, note that $\tilde w(\alpha) \equiv w(\alpha) \Def (1-\alpha)w$, and 
\beq\label{eq-OAlpha}
\ba{rcl}
\rho(w(\alpha)) & \refEQ{eq-OLow} & {1 - \alpha \over n+1}(v_0 - (1-\alpha) \| v \|^2 ) \; = \; (1-\alpha) \rho(w) + {\alpha(1-\alpha) \over n+1} \| v \|^2.
\ea
\eeq
At the same time, we have
\beq\label{eq-First}
\ba{rcl}
\check r(\tilde z(\alpha)) & = & \tilde x(\alpha) \tilde s(\alpha) - \tilde v^2(\alpha) \; = \; (x + \alpha \tilde \Delta^x) (s + \alpha \tilde \Delta^s) - (1-\alpha)^2 v^2\\
\\
& = & x s + \alpha^2 \tilde \Delta^x\tilde \Delta^s  - (1-\alpha)^2 v^2 + \alpha \Big[ \frac{\| v \|^2 }{ n+1}\check e - v^2 - xs\Big]\\
\\
& \refEQ{eq-OAlpha} & \rho(w(\alpha))\check e + (1-\alpha) [xs - v^2 - \rho(w) \check e] + \alpha^2 \Big[  \frac{\| v \|^2 }{ n+1}\check e + \tilde \Delta^x\tilde \Delta^s - v^2 \Big].
\ea
\eeq
Finally,
$$
\ba{rcl}
r^{(0)}(\tilde z(\alpha)) & = & (1- \alpha) v^{(0)} - \la \tilde s(\alpha), \tilde x(\alpha)\ra \; = \; (1-\alpha)(v^{(0)} - \la s, x \ra) + {\alpha\over n+1} \| v \|^2\\
\\
& \refEQ{eq-OAlpha} & \rho(w(\alpha)) + (1-\alpha) (v^{(0)} - \la s, x \ra - \rho(w)) + {\alpha^2 \over n+1} \| v \|^2.
\ea
$$
Now we can see the main advantage of the direction $\tilde d$: the initial residual $r(z) - \rho(w) e$ is eliminated {\em automatically} by big steps. This opens a possibility of making large steps.

Thus, we have proved the following representation:
\beq\label{eq-RepFM}
\ba{rcl}
\Psi(\tilde z(\alpha)) & = & - \sum\limits_{i=0}^n \ln \left(1 + {1 \over \rho(w(\alpha))} B^{(i)}(\alpha) \right), \\
\\
B(\alpha) & \Def & (1-\alpha)(r(z) - \rho(w) e) + \alpha^2 \tilde g,\\
\\
\tilde g^{(0)} & = & \frac{\| v \|^2 }{ n+1}, \quad \tilde g^{(i)} \; = \; \frac{\| v \|^2 }{ n+1} + [\tilde \Delta^x\tilde \Delta^s- v^2]^{(i)}, \quad i = 1, \dots, n.
\ea
\eeq
\BL\label{lm-LStep}
Let point $z=(u,w) \in {\cal F}$ satisfy the centering condition $\delta(z) \leq \beta$. If the parameter $\alpha$ is chosen in accordance to the rules of Step b) of method (\ref{met-AutoPT}), then
\beq\label{eq-AStep}
\ba{rcl}
\half r & \leq & {\alpha^2 \over (1 - \alpha)\rho(w) } \| \tilde g \|.
\ea
\eeq
\EL
\proof
Note that $\la e, B(\alpha) \ra \refEQ{eq-Sum} 0$.
Assuming that ${1 \over \rho(w(\alpha))} \| B(\alpha)\| < r$, we get
$$
\ba{rcl}
\omega_*(r) & = & A_{\psi} \; \refEQ{met-AutoPT} \; \Psi(\tilde z(\alpha)) \; \refEQ{eq-RepFM} \; - \sum\limits_{i=0}^n \ln \left(1 + {1 \over \rho(w(\alpha))} B^{(i)}(\alpha) \right) \\
& \leq & \omega_* \left( {1 \over\rho(w(\alpha))} \| B(\alpha) \| \right) \; < \; \omega_*(r),
\ea
$$
which is impossible. Therefore, ${1 \over \rho(w(\alpha))} \| B(\alpha)\| \geq r$.

Since $\delta(z) \leq \beta$, by the second inequality in (\ref{eq-Rho}), we have
\beq\label{eq-Rho1}
\ba{rcl}
\chi_0(z) & = & {1 \over \rho(w)} \| r(z) - \rho(w) e \| \; \leq \; {\beta \over 1- \beta}.
\ea
\eeq
Hence 
$$
\ba{rcl}
r & \leq & {1 \over \rho(w(\alpha))} \| B(\alpha)\| \; \refLE{eq-OAlpha} \; { 1 \over (1-\alpha) \rho(w)} \Big[ (1-\alpha) \| r(z) - \rho(w) e \| + \alpha^2 \| \tilde g \| \Big]\\
\\
& \refLE{eq-Rho1} & {\beta \over 1 - \beta} + {\alpha^2 \over (1-\alpha) \rho(w)} \| \tilde g \| \; = \; \half r + {\alpha^2 \over (1-\alpha) \rho(w)} \| \tilde g \|. \QR
\ea
$$

Inequality (\ref{eq-AStep}) is our main tool in the convergence analysis of method (\ref{met-AutoPT}). For the local convergence, we use its simplified version:
\beq\label{eq-ALoc}
\ba{rcl}
1 - \alpha & \leq & {2 \over \rho(w) r} \| \tilde g \|.
\ea
\eeq
Thus, we need to find an upper bound for $\| \tilde g \|$.
Note that the results of Section \ref{sc-Loc} are valid for any right-hand side $d$. Hence, 
$$
\ba{rcl}
\| \tilde \Delta^x \tilde \Delta^s \| & \refLE{eq-SumBND} &  {2^{3/2}\kappa \over \pi_*} \| \tilde d \|^2.
\ea
$$
Thus, for $v_+ = (0,v) \in \R^{n+1}$, we have
$$
\ba{rcl}
\| \tilde g \| & \refLE{eq-RepFM} & \Big\| {\| v \|^2\over n+1} e - v_+^2 \Big\| + \| \tilde \Delta^x \tilde \Delta^s \| \; \refLE{eq-DV4} \; \| v \|^2 +  {2^{3/2}\kappa \over \pi_*} \| \tilde d \|^2.
\ea
$$
At the same time,
\beq\label{eq-TLD}
\ba{rcl}
\half \| \tilde d \|^2 & \refLE{eq-DTSum} & \| d_a \|^2 + \| d_c \|^2 \; \refLE{eq-DA} \; {(v^{(0)})^2 \over n+1} +  {4 v^{(0)} \| v \|^2 \over (n+1)^2}  + 4\| v \|^4 + \| r(z) - \rho(w) e \|^2\\
\\
& \refLE{eq-Rho1} & {(v^{(0)})^2 \over n+1} +  {4 v^{(0)} \| v \|^2 \over (n+1)^2}  + 4\| v \|^4 + \rho^2(w) {\beta^2 \over (1-\beta)^2}  \\
\\
& \leq & \left({1 \over n+1} + {4 \over (n+1)^2} + 4 + {\beta^2 \over (n+1)^2 (1-\beta)^2} \right) (v^{(0)})^2.
\ea
\eeq
Putting all inequalities together, we get
$$
\ba{rcl}
1 - \alpha & \leq & {2 \over \rho(w) r} \Big[ \| v \|^2 +  {2^{5/2}\kappa \over \pi_*} \left( {(v^{(0)})^2 \over n+1} +  {4 v^{(0)} \| v \|^2 \over (n+1)^2}  + 4\| v \|^4 + \rho^2(w) {\beta^2 \over (1-\beta)^2} \right) \Big]\\
\\
& = & {2 \over \rho(w) r} \Big[ \| v \|^2 + { 2^{5/2}\kappa \over \pi_*} \left( {(v^{(0)})^2 \over n+1} +  {4 v^{(0)} \| v \|^2 \over (n+1)^2}  + 4\| v \|^4  \right) \Big] + { 2^{3/2} \kappa \over  \pi_*} r \rho(w).
\ea 
$$

Coming back to the whole iteration process, we denote $f_k = w_k^{(0)} = \tau_k \bar v_0$, where $(\bar v_0, \bar v) \equiv w_0$. Then, as in the relations (\ref{eq-BB1}), we get
$$
\ba{rcl}
1 - \alpha_k & \leq & {2 (n+1) \tau_k \over (\bar v_0 - \tau_k \| \bar v \|^2)r}\Big[ \| \bar v \|^2 +  {2^{5/2} \kappa \over \pi_*} \left( {\bar v_0^2 \over n+1} +   {4 \bar v_0 \| \bar v \|^2 \tau_k \over (n+1)^2}  + 4 \tau_k^2 \| \bar v \|^4  \right) \Big] + { 2^{3/2}\kappa r \over  \pi_*} \tau_k {\bar v_0 - \tau_k \| \bar v \|^2 \over n+1}\\
\\
& \leq & {2 (n+1) f_k \over (f_0 - f_k)r}\Big[ 1 +  {2^{5/2}\kappa \over \pi_*} \left( {f_0 \over n+1} +  {4 f_k \over (n+1)^2}  + 4 {f_k^2 \over f_0} \right) \Big] + { 2^{3/2}\kappa r f_k \over  \pi_*(n+1)}\\
\\
& \leq &  {2 (n+1) f_k \over (f_0 - f_k)r}\Big[ 1 + 22 \sqrt{2}{\kappa \over \pi_*} f_0 \Big] + { 2^{3/2} \kappa r f_k \over  \pi_*(n+1)}.
\ea
$$

As for Theorem \ref{th-Quad}, we assume that the starting point $u_0$ satisfies condition $\la s_0, x_0 \ra \leq {\pi_* \over 4 \kappa}$. Thus, we have proved the following statement.
\BT\label{th-Quad2} Let $f_0 \leq {\pi_* \over 4 \kappa}$. Then, for the method (\ref{met-AutoPT}), we have
\beq\label{eq-2Rate}
\ba{rcl}
f_{k+1} & \leq & \Big[ {18 (n+1)  \over (f_0 - f_k)r} + { r \sqrt{2} \over 2 (n+1)f_0} \Big] f^2_k.
\ea
\eeq
\ET

Note that now we can keep the proximity level $\beta$ constant.
In the remaining part of this sections, using  the inequality (\ref{eq-AStep}), we justify the global complexity bound of the method (\ref{met-AutoPT}). For that, we need to find another upper bound for $\| \tilde g \|$. Since
$$
\ba{rcl}
\| \tilde g \| & \refLE{eq-RepFM} & \Big\| {\| v \|^2 \over n+1} e - v_+^2 \Big\| + \| \tilde \Delta^x \tilde \Delta^s \| \; \refLE{eq-DV4} \; \| v \|^2 + \| \tilde \Delta^x \tilde \Delta^s \|,
\ea
$$
we need to estimate the last term. Note that
\beq\label{eq-TL1}
\ba{rcl}
2 \| \tilde \Delta_x \tilde \Delta_s \| & \leq & 2 \sum\limits_{i=1}^n \Big| [\tilde \Delta^x \tilde \Delta^s]^{(i)} \Big| \; \leq \; \la S X^{-1} \tilde \Delta^x, \tilde \Delta^x \ra + \la X S^{-1} \tilde \Delta^s, \tilde \Delta^s \ra.
\ea
\eeq

Denote $\ua = {v_0 \over \| v \|^2} > 1$. 
\BL\label{lm-TL2}
We have the following bound:
\beq\label{eq-TL2}
\ba{rcl}
\la S X^{-1} \tilde \Delta^x, \tilde \Delta^x \ra + \la X S^{-1} \tilde \Delta^s, \tilde \Delta^s \ra & \leq & n_r \left( {\ua \over \ua - 1} \right)^2  \rho(w),
\ea
\eeq
where $n_r = {25 \over 6} + {n \over 1-\beta} = {25 \over 6} + \left(1 + {r \over 2} \right)n$.
\EL
\proof
Indeed, in view of equality $\la \tilde \Delta^s, \tilde \Delta^x \ra = 0$, we have
$$
\ba{rcl}
\la S X^{-1} \tilde \Delta^x, \tilde \Delta^x \ra + \la X S^{-1} \tilde \Delta^s, \tilde \Delta^s \ra & = & \|  X^{1/2} S^{-1/2} \tilde \Delta^s + S^{1/2} X^{-1/2} \tilde \Delta^x) \|^2 \\
\\
& \refEQ{eq-LSys} & \| X^{-1/2} S^{-1/2} \tilde d \|^2.
\ea
$$
Since $\tilde d = {v^{(0)} \over n+1} \check e - x(w)s(w) - xs = {v^{(0)} \over n+1} \check e - 2 xs + \check r(z) - \rho(w) \check e$, we have
$$
\ba{rl}
& \half \| X^{-1/2} S^{-1/2}\tilde d \|^2 \\
\\
\leq & \| X^{-1/2} S^{-1/2} ({v^{(0)} \over n+1} \check e - 2 xs) \|^2 + \| X^{-1/2}S^{-1/2} (\check r(z) - \rho(w) \check e)\|^2\\
\\
\leq & {(v^{(0)})^2 \over (n+1)^2} \sum\limits_{i=1}^n {1 \over x^{(i)} s^{(i)}} - 4 {n \over n+1} v^{(0)} + 4\la s, x \ra + \chi_1^2(z) \rho(w) \\
\\
\stackrel{(\ref{eq-PropXS}),(\ref{eq-Rho})}{\leq} & 4 {v^{(0)} \over n+1} + {n (v^{(0)})^2 \over (n+1)^2 (1-\beta) \rho(w)} + {\beta^2 \over 1-\beta} \rho(w) \\
\\
 = & \left( {4 v^{(0)} \over v^{(0)} - \| v \|^2} + {n (v^{(0)})^2 \over (1-\beta)(v^{(0)} - \| v \|^2)^2} + {\beta^2 \over 1-\beta} \right) \rho(w)\\
 \\
\leq & \left({1 \over 6} + {4 \ua \over \ua - 1} + {n \over 1 - \beta} \left( {\ua \over \ua - 1} \right)^2 \right) \bar \rho(w) \;
\leq \; n_r \left( {\hat \alpha \over \hat \alpha - 1} \right)^2  \rho(w). \QF
\ea
$$

Thus, we conclude that
$$
\ba{rcl}
{1 \over \rho(w)} \| \tilde g \| & \leq & {\| v \|^2 \over \rho(w)}  + \half n_r \left( {\ua \over \ua - 1} \right)^2 \; = \; {n+1 \over \ua - 1} + \half n_r \left( {\ua \over \ua - 1} \right)^2\; \leq \; ({1 + n \over 4} + \half n_r)\left( {\ua \over \hat \ua - 1} \right)^2.
\ea
$$
Hence, denoting $\tilde n_r = {n+1 \over 2} + n_r$, we conclude that $\alpha$ satisfies the following inequality
\beq\label{eq-AG}
\ba{rcl}
r \tilde n_r^{-1} & \leq & {\alpha^2  \over 1 - \alpha}  \left( {\ua \over \ua - 1} \right)^2.
\ea
\eeq
\BL\label{lm-AG}
Let $\alpha > 0$ satisfy inequality (\ref{eq-AG}). Then $\alpha \geq \gamma {\hat \alpha - 1 \over \hat \alpha}$ with $\gamma = { 1 \over 1 + \sqrt{ \tilde n_r/r }}$.
\EL
\proof
Indeed, from inequality (\ref{eq-AG}), we have $r \tilde n_r^{-1}\leq \left( {\ua \alpha \over (1-\alpha)(\ua - 1)} \right)^2$.
Hence, ${\alpha \over 1 - \alpha} \geq {\ua - 1 \over \ua} \sqrt{r \tilde n_r^{-1}} $. Therefore,
$\alpha \geq { {\ua - 1 \over \ua} \sqrt{r \tilde n_r^{-1}} \over 1 + {\ua - 1 \over \ua} \sqrt{r \tilde n_r^{-1}}} \; = \; { \ua - 1   \over \ua \sqrt{ \tilde n_r /r }+ \ua - 1} \; \geq \; { \ua - 1   \over \ua \left(1 + \sqrt{ \tilde n_r/r } \right)}$. 
\qed

Thus, in view of Lemma 5.7 in \cite{AFF}, method (\ref{met-AutoPT}) has the following rate of convergence:
\beq\label{eq-Rate}
\ba{rcl}
\mu^*(w_{k+1}) & \leq & {1 \over 1 + \gamma} \mu^*(w_k),
\ea
\eeq
where $\mu^*(w) = {(v^{(0)})^2 \over v^{(0)} - \| v \|^2} \geq v^{(0)}$ and $\gamma = {1 \over 1 + \sqrt{ \tilde n_r/r }}$. Note that 
$$
\ba{rcl}
{1 \over r} \tilde n_r \; \approx \; {1 \over r}({n \over 2} + {n \over 1 - \beta}) & = & n \cdot {1-\beta \over 2 \beta} \cdot {3 - \beta \over 2(1-\beta)} \; = \; {3-\beta \over 4 \beta} n. 
\ea
$$
For the choice $r = {6 \over 7}$, we have $A_{\psi} \approx 1.09$ and $\beta = 0.3$. Therefore, 
${1 \over r} \tilde n_r \approx {9 \over 4}n$, and 
\beq\label{eq-Gamma1}
\ba{rcl}
\gamma & \approx & {2 \over 3 \sqrt{n}}.
\ea
\eeq

\section{Second-order prediction}\label{sc-Second}
\SetEQ

Let us include in Step b) of method (\ref{met-AutoPT}) a second-order prediction. 
\beq\label{met-PT2}
\ba{|l|}
\hline \\
\hspace{8ex} \mbox{\bf $2^{\mbox{nd}}$-order Prediction for PTFM-Method (PTFM2)} \quad \\
\\
\hline \\
\mbox{{\bf Initialization.} Choose  $r \in (0,1)$, $A_{\psi} = \omega_*(r)$, $u_0 \in {\cal F}_0$, and $w_0 \refEQ{eq-Start} w_*(u_0)$.}\\
\mbox{Define the acceptance level $\beta = {r \over 2+r} < {1 \over 3}$.}\\
\\
\mbox{\bf $k$th iteration ($k \geq 0$).}\\
\\
\ba{rl}
\mbox{{\bf a)}} & \mbox{Compute $r(z_k)$ and $\Sigma_k^{-1} = \left[A X_k S_k^{-1} A^T\right]^{-1}$. }\\
\\
\mbox{{\bf b)}} & \mbox{If $\delta(z_k) \leq \beta$, then do \hspace{2ex} \fbox{\sc Predictor Step}}\\
& \bullet \; \mbox{Set $\tilde d_k = \frac{\| v_k \|^2 }{ n+1} \check e - v_k^2 - x_k s_k$ and compute $\tilde \Delta_k = \Delta_k(\tilde d_k)$.}\\
& \bullet \; \mbox{Set $\widehat d_k = v_k^2 - \frac{\| v_k \|^2 }{ n+1} \check e - \tilde \Delta^x_k \tilde \Delta^s_k$ and compute $\widehat \Delta_k = \Delta_k(\widehat d_k)$.}\\
& \bullet \; \mbox{Define function  $\widehat \psi_k(\alpha) = \Psi(u_k+\alpha \tilde \Delta_k + \alpha^2 \wh \Delta_k, (1-\alpha)w_k)$.}\\
& \bullet \; \mbox{Find $\alpha_k$ as an approximate solution of equation $\widehat \psi_k(\alpha) = A_{\psi}$.}\\
& \bullet \; \mbox{Define $u_{k+1} =u_k + \alpha_k \tilde \Delta_k + \alpha_k^2 \wh \Delta_k$ and $w_{k+1} = (1-\alpha_k)w_k$.}\\
\\
\mbox{\bf c)} & \mbox{Otherwise, do \hspace{2ex} \fbox{\sc Corrector Step}}\\
& \bullet \; \mbox{Define $d_k = \rho(w_k) \check e - \check r(z_k)$. Compute $\Delta_k = \Delta_k(d_k)$.}\\
& \bullet \; \mbox{Define function  $f_k(\alpha) = F(u_k+\alpha \Delta_k, w_k)$.}\\
& \bullet \; \mbox{Find $\alpha_k$ as an approximate minimum of $f_k(\alpha)$ in $\alpha \geq 0$.}\\
& \bullet \; \mbox{Define $u_{k+1} =u_k + \alpha_k \Delta_k$ and $w_{k+1} = w_k$.}\\
\\
\mbox{\bf d)} & \mbox{If $w_k^{(0)} \leq \epsilon$ and $\delta(z_k) \leq \beta$, then \fbox{\sc Stop}}
\ea\\
\\
\hline
\ea
\eeq

Let us analyze the predictor Step b) of method (\ref{met-PT2}). In our reasoning, for the sake of notation, we omit index $k$. Denote $\wh z(\alpha) = z + \alpha (\tilde \Delta,-w) + \alpha^2 (\wh \Delta,0)$. Then,
$$
\ba{c}
\check r(\wh z(\alpha)) = (x + \alpha \tilde \Delta^x+\alpha^2 \wh \Delta^x)(s + \alpha \tilde \Delta^s+\alpha^2 \wh \Delta^s) - (1-\alpha)^2 v^2\\
\\
= (x + \alpha \tilde \Delta^x)(s + \alpha \tilde \Delta^s) + \alpha^2 [\wh \Delta^x (s + \alpha \tilde \Delta^s)+\wh \Delta^s (x + \alpha \tilde \Delta^x)] + \alpha^4 \wh \Delta^x\wh \Delta^s- (1-\alpha)^2 v^2\\
\\
\refEQ{eq-First} \rho(w(\alpha)) \check e + (1-\alpha)[\check r(z)- \rho(w) \check e] + \alpha^3 [ \wh \Delta^x \tilde \Delta^s + \wh \Delta^s \tilde \Delta^x ] + \alpha^4 \wh \Delta^x\wh \Delta^s.
\ea
$$
Similarly,
$$
\ba{rcl}
r^{(0)}(\wh z(\alpha)) & = & (1-\alpha)v^{(0)} - \la s + \alpha \tilde \Delta^s + \alpha^2 \wh \Delta^s, x + \alpha \tilde \Delta^x + \alpha^2 \wh \Delta^x \ra\\
\\
& = & (1-\alpha)v^{(0)} - \la s,x \ra - \alpha [ \la s, \tilde \Delta^x \ra + \la \tilde \Delta^s, x \ra] - \alpha^2 [ \la s, \wh \Delta^x \ra + \la \wh \Delta^s, x \ra]\\
\\
& = & (1-\alpha)v^{(0)} - \la s,x \ra - \alpha [ -{\| v \|^2 \over n+1} - \la s,x \ra] - \alpha^2 {1 \over n+1} \| v \|^2\\
\\
& = & (1-\alpha)r^{(0)}(z) + \alpha(1-\alpha) {1 \over n+1} \| v \|^2\\
\\
& \refEQ{eq-OAlpha} & \rho(w(\alpha)) + (1-\alpha)(r^{(0)}(z) - \rho(w)).
\ea
$$
Thus, we have proved the following representation:
\beq\label{eq-Rep2}
\ba{rcl}
\Psi(\wh z(\alpha)) & = & - \sum\limits_{i=0}^n \ln \Big( 1 + {1 \over \rho(w(\alpha))} C^{(i)}(\alpha) \Big), \\
\\
C(\alpha) & = & (1-\alpha) (r(z) - \rho(w) e) + \alpha^3 g_1 + \alpha^4 g_2,\\
\\
g_1 & = & \left( \ba{c} 0 \\\wh \Delta^x \tilde \Delta^s + \wh \Delta^s \tilde \Delta^x \ea \right), \quad g_2 \; = \; \left( \ba{c} 0 \\\wh \Delta^x \wh \Delta^s \ea \right).
\ea
\eeq
Note that $\la e, C(\alpha) \ra \equiv 0$. Hence, assuming that ${1 \over \rho(w(\alpha))}\|  C(\alpha) \| < r$, we get
$$
\ba{rcl}
\omega_*(r) & > & \Psi(\wh z(\alpha)) \; = \; A_{\psi} \; = \; \omega_*(r),
\ea
$$
which is impossible. Hence, we conclude that
$$
\ba{rcl} 
r & \leq & {1 \over \rho(w(\alpha))}\|  C(\alpha) \| \; \leq \; {1 \over (1-\alpha) \rho(w)} \Big[ (1-\alpha) \| r(z) - \rho(w) e\| + \alpha^3 \| g_1 \| + \alpha^4 \| g_2 \| \Big]\\
\\
& = & \chi_0(z) + {\alpha^3 \over (1-\alpha) \rho(w)} \Big[ \| g_1 \| + \alpha \| g_2 \| \Big] \; \refLE{eq-Rho} \; {\beta \over 1 - \beta} + {\alpha^3 \over (1-\alpha) \rho(w)} \Big[ \| g_1 \| + \alpha \| g_2 \| \Big].
\ea
$$
Thus, in view of the choice of $\beta$ in (\ref{met-PT2}), we have
\beq\label{eq-LAL3}
\ba{rcl}
{r \over 2} & \leq & {\alpha^3 \over (1-\alpha) \rho(w)} \Big[ \| g_1 \| + \alpha \| g_2 \| \Big].
\ea
\eeq
For estimating the local convergence, we need a relaxed version of this inequality:
\beq\label{eq-1AL}
\ba{rcl}
1 - \alpha & \leq & {2 \over \rho(w) r} \Big[ \| g_1 \| + \| g_2 \| \Big].
\ea
\eeq

Let us estimate now the norms of vectors $g_1$ and $g_2$, assuming that the condition (\ref{eq-Neib}) is satisfied. Note that
$$
\ba{rcl}
\| g_1 \| + \| g_2 \| & \refLE{eq-Prod} &  \| \wh \Delta^x \| \cdot \| \tilde \Delta^s \| + \| \wh \Delta^s \| \cdot \| \tilde \Delta^x \| + \| \wh \Delta^s \| \cdot \| \wh \Delta^x \| \\
\\
& \refLE{eq-TBound} & {5 \over 2 \pi_*}(1+\kappa^2) \Big[ 2 \| \hat d \| \, \| \tilde d \| + \| \hat d \|^2 \Big].
\ea
$$

At the same time,
$$
\ba{rcl}
\half \| \tilde d \|^2 & \refLE{eq-TLD} & \left({11 \over 2} + {r^2 \over 16} \right) (v^{(0)})^2,\\
\\
\| \hat d \| & \leq & \Big\| v^2 - {\| v \|^2 \over n+1} \check e \Big\| + \| \tilde \Delta^x \tilde \Delta^s \| \; \stackrel{(\ref{eq-DV4}),(\ref{eq-SumBND})}{\leq} \; \| v \|^2 + {2^{3/2}\kappa \over \pi_*} \| \tilde d \|^2\\
\\
& \leq & \| v \|^2 + {2^{3/2}\kappa \over \pi_*} \left(11 + {r^2 \over 8} \right) (v^{(0)})^2.
\ea
$$
Denoting $c_1 = \sqrt{11 + {r^2 \over 8}}$ and $c_2 =  {2^{3/2}\kappa \over \pi_*} c_1^2$, we have
$$
\ba{rcl}
\| \tilde d \| & \leq & c_1 v^{(0)}, \quad \| \hat d \| \; \leq \; \| v \|^2 + c_2 (v^{(0)})^2.
\ea
$$
Denoting now $f_k = w_k^{(0)} = \tau_k f_0$, we get
$$
\ba{rcl}
f_{k+1} & \leq & (1-\alpha_k) f_k \; \refLE{eq-1AL} \; {2 (n+1) f_k \over r \tau_k(f_0 - f_k)} \cdot {5(1+\kappa^2) \over 2 \pi_*} \Big[ 2 c_1 \tau_k^3 f_0 (f_0 + c_2 f_0^2) + \tau_k^4 (f_0 + c_2 f_0^2)^2 \Big]\\
\\
& = & {5(n+1)(1+\kappa^2) f^3_k \over \pi_* r (f_0 - f_k)}  \Big[ 2 c_1 (1 + c_2 f_0) + \tau_k (1 + c_2 f_0)^2 \Big].
\ea
$$
Since our estimates are valid only for $f_0 \leq {\pi_* \over 4 \kappa}$, we conclude that
\beq\label{eq-Rate3}
\ba{rcl}
f_{k+1} & \leq & {5(n+1) (1+\kappa^2) c_3 \over  \pi_* r (f_0 - f_k)} \, f^3_k \; \leq \; (n+1)\left({1 \over \kappa}+\kappa\right) {5c_3 f_k^3 \over 4f_0 (f_0 - f_k)} \, \\
\\
 c_3 & \Def & \left(1 + {c_1^2 \over \sqrt{2}} \right)
\left(1+2c_1+ {c_1^2 \over \sqrt{2}}\right).
\ea
\eeq

Let us prove now the polynomial-time complexity of method (\ref{met-PT2}). First of all, we need to justify the following bound. 
\BL\label{lm-Hat}
Under condition of Step b) in method (\ref{met-PT2}), we have
\beq\label{eq-Hat}
\ba{rcl}
\la S X^{-1} \wh \Delta^x, \wh \Delta^x \ra + \la XS^{-1} \wh \Delta^s, \wh \Delta^s \ra & \leq & \Big[ 
{2(n+1)\ua \over (\ua - 1)^2} + {n_r^2 \over 2} 
\left({\ua \over \ua - 1} \right)^4 \Big] {\rho(w) \over 1 - \beta} \\
\\
& \leq & \hat n_r^2 \left({\hat \alpha \over \hat \alpha - 1} \right)^4 {\rho(w) \over 1 - \beta},
\ea
\eeq
where $\hat n^2_r = {16 \over 27} (n+1) + \half n_r^2$.
\EL
\proof
Note that $\la \wh \Delta^s, \wh \Delta^x \ra = 0$. Therefore, 
$$
\ba{rcl}
\la S X^{-1} \wh \Delta^x, \wh \Delta^x \ra + \la XS^{-1} \wh \Delta^s, \wh \Delta^s \ra & = & \| S^{1/2} X^{-1/2} \wh \Delta^x + X^{1/2} S^{-1/2} \wh \Delta^s \|^2 \\
\\
& = & \| (SX)^{-1/2} \hat d \|^2. 
\ea
$$
Since $\hat d = v^2 - {\| v \|^2 \over n+1} \check e - \tilde \Delta^x \tilde \Delta^s$, we have
$$
\ba{rcl}
\half \| (SX)^{-1/2} \hat d \|^2 & \leq & \Big\| (SX)^{-1/2}(v^2 - {\| v \|^2 \over n+1} \check e) \Big\|^2 + \| (SX)^{-1/2}\tilde \Delta^x \tilde \Delta^s \|^2.
\ea
$$
Note that $xs\refGE{eq-PropXS} (1-\beta)x(w)s(w)$, with $x(w) s(w) = v^2 + \rho(w) \check e \geq \rho(w) \check e$. Hence,
$$
\ba{rl}
& \Big\| (SX)^{-1/2}(v^2 - {\| v \|^2 \over n+1} \check e) \Big\|^2 \; = \; \Big\| (SX)^{-1/2}(x(w)s(w) - {v^{(0)} \over n+1} \check e) \Big\|^2\\
\\
\leq & {1 \over 1-\beta} \Big\| (S(w)X(w))^{-1/2}(x(w)s(w) - {v^{(0)} \over n+1} \check e) \Big\|^2 \\
\\
 = & {1 \over 1 - \beta} \Big[ \la s(w), x(w) \ra - 2{n \over n+1} v^{(0)} + \left({v^{(0)} \over n+1} \right)^2 \sum\limits_{i=1}^n {1 \over x^{(i)}(w) s^{(i)}(w)} \Big]\\
 \\
 \leq & {1 \over 1 - \beta} \Big[ \| v \|^2 + n \rho(w) - 2{n \over n+1} v^{(0)} + \left({v^{(0)} \over n+1} \right)^2 {n \over \rho(w)} \Big]\\
 \\
 = & {1 \over 1 - \beta} \Big[ {1 \over n+1}\| v \|^2  - {n \over n+1} v^{(0)} + {nv_0^2 \over (n+1)(v^{(0)} - \| v \|^2)} \Big] \; = \; {1 \over 1 - \beta} \Big[ {1 \over n+1}\| v \|^2  + {nv^{(0)} \| v \|^2 \over (n+1)(v^{(0)} - \| v \|^2)} \Big]\\
 \\
 = & {1 \over 1 - \beta} \Big[ {\rho(w) \over \ua - 1} + {n \over n+1} \cdot {\ua \over \ua - 1} \cdot {(n+1)\rho(w) \over \ua - 1} \Big] \; = \; {\rho(w) \over 1 - \beta} \Big[{1 \over \ua - 1} +  {\ua n \over (\ua-1)^2}\Big] \; \leq \;  {\rho(w) \over 1 - \beta} \cdot{\ua (n+1) \over (\ua-1)^2}.
\ea
$$

For the second term, we have
$$
\ba{rcl}
\| (SX)^{-1/2}\tilde \Delta^x \tilde \Delta^s \|^2 & \leq & {1 \over (1-\beta) \rho(w)} \| \tilde \Delta^x \tilde \Delta^s \|^2 \\
\\
& \refLE{eq-TL1} & {1 \over (1-\beta) \rho(w)} \Big[ \half \la SX^{-1} \tilde \Delta^x, \tilde \Delta^x \ra + \half \la XS^{-1} \tilde \Delta^s, \tilde \Delta^s \ra \Big]^2 \\
\\
& \refLE{eq-TL2} & {1 \over (1-\beta) \rho(w)} \Big[ \half n_r \left({\ua \over \ua - 1} \right)^2 \rho(w)\Big]^2 \; = \; {n_r^2 \over 4(1-\beta)} \left({\ua \over \ua - 1} \right)^4 \rho(w).
\ea
$$

Thus, it remains to combine the bounds for two terms.
\qed

Now we can estimate the norms of the vectors $g_1$ and $g_2$.
$$
\ba{rcl}
\| g_2 \| & = & \| \wh \Delta^x \wh \Delta^s \| \; \leq \; \sum\limits_{i=1}^n \Big| [\wh \Delta^x \wh \Delta^s]^{(i)} \Big| \; \leq \; 
\half \la S X^{-1} \wh \Delta^x, \wh \Delta^x \ra + \half \la XS^{-1} \wh \Delta^s, \wh \Delta^s \ra \\
& \refLE{eq-Hat} & {\hat n_r^2 \over 2(1-\beta)} \left({\ua \over \ua - 1} \right)^4 \rho(w).
\ea
$$
For the first vector $g_1$, let us choose a scaling coefficient $\tau > 0$. Then
$$
\ba{rl}
\| g_1 \| \leq & \| \wh \Delta^s \tilde \Delta^x \| + \| \wh \Delta^x \tilde \Delta^s \| \; \leq \; \sum\limits_{i=1}^n \Big[\Big| [\wh \Delta^s \tilde \Delta^x]^{(i)} \Big| + \Big| [\wh \Delta^x \tilde \Delta^s]^{(i)} \Big| \Big]\\
\\
\leq & \half \sum\limits_{i=1}^n \Big[\tau {x^{(i)} \over s^{(i)}} (\wh \Delta^s \wh \Delta^s)^{(i)} + {s^{(i)} \over \tau x^{(i)}} (\tilde \Delta^x \tilde \Delta^x)^{(i)} + \tau {s^{(i)} \over x^{(i)}} (\wh \Delta^x \wh \Delta^x)^{(i)} +  {x^{(i)} \over \tau s^{(i)}}(\tilde \Delta^s \tilde \Delta^s)^{(i)} \Big]\\
\\
= & {\tau \over 2} \left[ \la XS^{-1} \wh \Delta^s, \wh \Delta^s \ra + \la SX^{-1} \wh \Delta^x, \wh \Delta^x \ra \right] + {1 \over 2 \tau} \left[ \la XS^{-1} \tilde \Delta^s, \tilde \Delta^s \ra + \la SX^{-1} \tilde \Delta^x, \tilde \Delta^x \ra \right]\\
\\
\stackrel{(\ref{eq-Hat}),(\ref{eq-TL2})}{\leq} &  {\tau \hat n_r^2 \over 2(1 - \beta)} \left({\hat \alpha \over \hat \alpha-1} \right)^4 \rho(w) + {n_r  \over 2 \tau} \left({\hat \alpha \over \hat \alpha-1} \right)^2 \rho(w).
\ea
$$
Minimizing the right-hand side of this inequality in $\tau$, we get the following bound:
$$
\ba{rcl}
\| g_1 \| & \leq & {\hat n_r n_r^{1/2} \over \sqrt{1 - \beta}} \left({\ua \over \ua-1} \right)^3 \rho(w) \; = \; {1 \over \sqrt{1 - \beta}} \left({\tilde n_r^{1/2}\ua \over \ua-1} \right)^3 \rho(w),
\ea
$$
where $\tilde n_r = \hat n^{2/3}_r n_r^{1/3}$.
Substituting these bounds in inequality (\ref{eq-LAL3}), we come to the following consequence:
\beq\label{eq-Step2}
\ba{rcl}
{r \over 2} & \leq & {1 \over \sqrt{1 - \beta}} \left({\bar n_r^{1/2}  \alpha \ua \over (1-\alpha)(\ua-1)} \right)^3 + {1 \over 2(1-\beta)} \left({ \bar n_r^{1/2} \alpha \hat \alpha \over (1-\alpha)(\hat \alpha - 1)} \right)^4,
\ea
\eeq
where $\bar n_r = \max\{ \hat n_r, n_r \}$.
Denoting $\tau = {1 \over 2\sqrt{1 - \beta}} \left({\bar n_r^{1/2}  \alpha \ua \over (1-\alpha)(\ua-1)} \right)$ and $\hat r = {r \over 16(1-\beta)}$, we get inequality $\hat r \leq \tau^3 + \tau^4$. Denote by $\tau_*$ the exact solution of the equation $\tau_*^3 + \tau_*^4 = \hat r$. Then $\tau \geq \tau_*$. Since $\tau_* \leq \hat r^{1/3}$, we have
$$
\ba{rcl}
\tau \; \geq \; \tau_* & = & {\hat r^{1/3} \over (1+\tau_*)^{1/3}} \; \geq \; {\hat r^{1/3} \over 1 + {1 \over 3} \hat r^{1/3}}.
\ea
$$
Thus, we come to the bound ${\alpha \over 1 - \alpha} \geq  
{2 \hat r^{1/3} \sqrt{1-\beta} \over 1 + {1 \over 3} \hat r^{1/3}} 
{\ua -1\over \ua}  \bar n_r^{-1/2} =  {\kappa_1(\ua -1)\over \kappa_2 \ua}  \bar n_r^{-1/2} $, where
$$
\ba{rcl}
\kappa_1 & = & \left({r \over 2} \sqrt{1-\beta} \right)^{1/3}, \quad \kappa_2 \; = \; 1 + {1 \over 3} \hat r^{1/3} \; = \; 1 + {1 \over 6} \left({r \over 2(1-\beta)}\right)^{1/3}.
\ea
$$
This bound can be rewritten as follows:
$$
\ba{rcl}
\alpha & \geq & {\kappa_1  (\ua -1) \over \kappa_1 (\ua -1) + \kappa_2 \ua \bar n_r^{1/2} } \geq {\kappa_1 (\ua - 1) \over (\kappa_1 + \kappa_2 \bar n_r^{1/2}) \ua}.
\ea
$$
Hence, the sequence of points generated by method (\ref{met-PT2}) satisfies inequality (\ref{eq-Rate}) with
\beq\label{eq-Rate2}
\ba{c}
\gamma \; = \; \gamma_2 \Def \; {\left({r \over 2} \sqrt{1-\beta} \right)^{1/3} \over \bar n_r^{1/2} \left( 1 + {1 \over 6} \left( {r \over 2(1-\beta)} \right)^{1/3} \right) +
\left({r \over 2} \sqrt{1-\beta} \right)^{1/3}},  \quad \bar n_r = \max\{ \hat n_r, n_r\}, \\
n_r \; = \; {25 \over 6} + {n \over 1-\beta} , \quad \hat n_r\;  = \; \sqrt{{16 \over 27}(n+1) + \half n_r^2} \; \approx\; {n \over (1-\beta)\sqrt{2}} \; < \; n_r.
\ea
\eeq
Thus, asymptotically, $\bar n_r = \hat n_r$ and the convergence rate of this scheme is defined by
$$
\ba{rcl}
\gamma_2 & \approx & 
{\left({r \over 2} \sqrt{1-\beta} \right)^{1/3} \over \hat n_r^{1/2} \left(1 + {1 \over 6} \left( {r \over 2(1-\beta)} \right)^{1/3} \right)}
 \; \approx \; {\beta^{1/3} (1-\beta) \over n^{1/2} \left((1-\beta)^{2/3} + {1 \over 6} \beta^{1/3} \right)}.
\ea
$$
For the recommended value $\beta =0.3$ (this is $r={6 \over 7}$), the coefficient $\gamma_2$ approaches ${0.52 \over n^{1/2}}$. It is slightly worse than the coefficient $\gamma_1 \refEQ{eq-Gamma1} {2 \over 3 n^{1/2}}$ for method (\ref{met-AutoPT}). However, the second-order scheme~(\ref{met-PT2}) has an advantage of faster local convergence. At the same time, the computational efforts required for one iteration in both methods are essentially the same.

\section{Finite termination}\label{sc-Fin}
\SetEQ

Local quadratic and cubic rates of convergence, presented in Sections \ref{sc-Auto} and \ref{sc-Second}, are so fast that for practical computations they are almost equivalent to finite termination of the corresponding schemes. It is interesting that, at the same time, the parabolic target-following methods can be endowed with a natural finite termination procedures, which need even less restrictive conditions than in Assumption \ref{sc-Loc}. This is the subject of the present section.

Our finite termination procedures are based on ordering of components of some {\em indicator vectors}, related to a particular parabolic target-following method. For $z = (x,s,y)$ in ${\cal F}_0$, we consider three different indicator vectors:
\begin{itemize}
\item Primal indicator vector $x$.
\item Dual indicator vector $s^{-1}$.
\item Primal-dual indicator vector $xs^{-1}$.
\end{itemize}
For a particular indicator vector $a \in \R^n_{++}$, denote by $\pi_a[\cdot]: [ 1: n] \to [ 1 : n ]$ the permutation function, representing the components of $a$ in a decreasing order:
$$
\ba{rcl}
a^{(\pi_a[i])} & \geq & a^{(\pi_a[i+1])}, \quad i = 1, \dots, n-1.
\ea
$$
With this function, we define the trial basis
$B_a = \{ k = \pi_a[i], \; i \in [1:m] \}$,
and compute the candidate optimal point $u^*_a = (x^*_a,s^*_a,y^*_a)$ in accordance to the following rules:
\beq\label{eq-Cand}
\ba{rcl}
x^*_{B_a} & = & A_{B_a}^{-1} b, \quad x^*_{N_a} \; = \; 0, \quad y_a^* \; = \; A_{B_a}^{-T} c_{B_a},\\
\\
s^*_{B_a} & = & 0, \quad s^*_{N_a} \; = \; c_{N_a} - A^T_{N_a} y_a^*,
\ea
\eeq
where $N_a = [1:n] \setminus B_a$. The test is successful if the matrix $A_{B_a}$ is non-degenerate and both vectors $x_a^* = x^*_{B_a} \bigcup x^*_{N_a}$, $s_a^* = s^*_{B_a} \bigcup s^*_{N_a}$ are non-negative.

Let us present the conditions, which guarantee that the point $u^*_a$ is indeed an optimal solution of the primal-dual problem (\ref{prob-LP}).
\BT\label{th-Fin}
Let problem (\ref{prob-LP}) has a unique optimal solution $u^* = (x^*,s^*, y^*)$ such that
\beq\label{eq-NonDG}
\ba{rcl}
x^* + s^* & > & 0.
\ea
\eeq
If the point $z = (u,w) \in {\cal F}$ satisfies the centering condition $\delta(z) \leq \beta$ with $\beta \in \left[0,{1 \over 3}\right)$, and
\beq\label{eq-PFin}
\ba{rcl}
\mu^*(w) = {(v^{(0)})^2 \over v^{(0)} - \| v \|^2} & < & (1-\beta) {\pi_* \over n+1},
\ea
\eeq
then the prediction $u^*_a$ formed by (\ref{eq-Cand}) with any of the indicator vectors $x$, $s^{-1}$, or $x s^{-1}$ is the optimal primal-dual solution of problem (\ref{prob-LP}).
\ET
\proof
Indeed, in view of the centering condition, we have
\beq\label{eq-XSLow1}
\ba{rcl}
x^{(i)} s^{(i)} & \refGE{eq-PropXS} & (1-\beta) x^{(i)}(w)s^{(i)}(w) \; \refGE{def-UW} \; (1-\beta) \rho(w), \quad i = 1, \dots, n.
\ea
\eeq
Denote by $B_*$ the set of positive components of $x^*$. Then, for any $i \in B_*$, we have
$$
\ba{rcl}
x^{(i)} & \geq & {1 - \beta \over s^{(i)}} \rho(w) \; \refGE{eq-Gap} \; {1 - \beta \over \la s, x \ra } \rho(w) x^*_{\min} \; \geq \; {1 - \beta \over (n+1) v^{(0)} } [v^{0)} - \| v \|^2] x^*_{\min}\\
\\
& \refGT{eq-PFin} & {v^{(0)} \over \pi_*} x^*_{\min} \; \geq \; {\la s, x \ra \over s^*_{\min}}.
\ea
$$
At the same time, for any $i \not\in B_*$, we have $x^{(i)} \refLE{eq-Gap} {\la s, x \ra \over s^*_{\min}}$. Hence, by ordering the components of vector $x$, we can detect the optimal basis $B_*$.

Similarly, for any $i \not\in B_*$ and $j \in B_*$, we have
$$
\ba{rcl}
s^{(i)} & \refGE{eq-XSLow1} & {1 - \beta \over x^{(i)}} \bar \omega(w) \; \refGE{eq-Gap} \; {1 - \beta \over \la s, x \ra } \bar \omega(w) s^*_{\min} \; \geq \; {1 - \beta \over (n+1) v^{(0)} } [v^{0)} - \| v \|^2] s^*_{\min}\\
\\
& \refGT{eq-PFin} & {v^{(0)} \over \pi_*} s^*_{\min} \; \geq \; {\la s, x \ra \over x^*_{\min}} \; \refGE{eq-Gap} \; s^{(j)}.
\ea
$$
Thus, by ordering components of vector $s$, we can detect the optimal basis $B_*$. The same is true for the vector $s^{-1}$.

Finally, since for both vectors $x$ and $s^{-1}$, the optimal basis corresponds to $m$ largest components, the same is true for the vector $xs^{-1}$.
\qed
\BC\label{cor-Complex}
Let problem (\ref{prob-LP}) satisfy the non-degeneracy assumption (\ref{eq-NonDG}). Then any of the methods (\ref{met-AffinePT}), (\ref{met-AutoPT}), or (\ref{met-PT2}), equipped with the termination procedure of Theorem~\ref{th-Fin}, can find its exact optimal solution in 
\beq\label{eq-Comp}
\ba{c}
O\left(\sqrt{n} \ln { (\la s_0,x_0 \ra + \sigma_0)^2 \over x^*_{\min} s^*_{\min} \sigma_0 } \right)
\ea
\eeq
iterations, where $u_0 = (x_0,s_0,y_0) \in {\cal F}_0$ is the starting point and $\sigma_0 = \min\limits_{1 \leq i \leq n} x^{(i)}_0 s^{(i)}_0$.
\EC
\proof
Indeed, $\mu^*(w_0) \refEQ{eq-Start} {(\la s_0,x_0 \ra + \sigma_0)^2 \over (n+1) \sigma_0}$. It remains to combine the condition (\ref{eq-PFin}) with the rate of convergence (\ref{eq-Rate}) and the lower bounds for the parameter $\gamma$ provided by Lemma \ref{lm-AG} and inequality (\ref{eq-Rate2}).
\qed

Since the main computational efforts at one iteration of the schemes (\ref{met-AutoPT}) and (\ref{met-PT2}) are spent for forming the matrix $\Sigma_k$, the optimality test (\ref{eq-Cand}) cannot not increase significantly the complexity of one iteration. However, for avoiding unnecessary computations at the first iterations, it may be reasonable to use the following activating conditions.
\BT\label{th-Awake}
Under conditions of Theorem \ref{th-Fin}, we have the following relations:
\beq\label{eq-AwakeX}
\ba{rcl}
\sum\limits_{i \in B_x} x^{(i)} & \geq & m^2 \sum\limits_{i \not\in B_x} x^{(i)},
\ea
\eeq
\beq\label{eq-AwakeS}
\ba{rcl}
\sum\limits_{i \not\in B_{1/s}} s^{(i)} & \geq & (n-m)^2 \sum\limits_{i \in B_{1/s}} s^{(i)},
\ea
\eeq
\beq\label{eq-AwakeXS}
\ba{rcl}
\sum\limits_{i \in B_{x/s}} {x^{(i)} \over s^{(i)}} & \geq & m^3 
\sum\limits_{i \not\in B_{x/s}} {x^{(i)} \over s^{(i)}}.
\ea
\eeq
\ET
\proof
In view of Theorem \ref{th-Fin}, we have $B_x = B_*$. Therefore,
$$
\ba{rcl}
\sum\limits_{i \in B_x} x^{(i)} & \refGE{eq-XSLow1} & (1-\beta) \rho(w) \sum\limits_{i \in B_x} {1 \over s^{(i)}} \\
& \refGE{eq-Gap} & (1-\beta) \rho(w) \min\limits_{s \geq 0} \left\{ \sum\limits_{i \in B_x} {1 \over s^{(i)}}: \; \sum\limits_{i \in B_x} s^{(i)} \leq {1 \over x^*_{\min}} \la s, x \ra \right\}\\
\\
& = & (1-\beta) \rho(w) {m^2 \over \la s, x \ra } x^*_{\min} \; \geq \; (1-\beta) {m^2 (v^{(0)} - \| v \|^2)\over (n+1) v^{(0)} } x^*_{\min} \; \refGE{eq-PFin} \; m^2 {v^{(0)} \over s^*_{\min}}.
\ea
$$
Since $\sum\limits_{i \not\in B_x} x^{(i)} \refLE{eq-Gap} {\la s, x \ra \over s^*_{\min}}$, we get (\ref{eq-AwakeX}).

Similarly, since $B_{1/s} = B_*$, we have
$$
\ba{c}
\sum\limits_{i \not\in B_{1/s}} s^{(i)} \; \refGE{eq-XSLow1} \; (1-\beta) \rho(w) \sum\limits_{i \not\in B_{1/s}} {1 \over x^{(i)}}\\
\refGE{eq-Gap} \; (1-\beta) \rho(w) \min\limits_{x \geq 0} \left\{ \sum\limits_{i \not\in B_{1/s}} {1 \over x^{(i)}}: \; \sum\limits_{i \not\in B_{1/s}} x^{(i)} \leq {1 \over s^*_{\min}} \la s, x \ra \right\}\\
\\
= \; (1-\beta) \rho(w) {(n-m)^2 \over \la s, x \ra } s^*_{\min} \; \geq \; (1-\beta) {(n-m)^2 (v^{(0)} - \| v \|^2)\over (n+1) v^{(0)} } s^*_{\min} \; \refGE{eq-PFin} \; (n-m)^2 {v^{(0)} \over x^*_{\min}}.
\ea
$$
Since $\sum\limits_{i \in B_{1/s}} s^{(i)} \refLE{eq-Gap} {\la s, x \ra \over x^*_{\min}}$, we get (\ref{eq-AwakeS}). Finally, we also have $B_{x/s} = B_*$. Therefore,
$$
\ba{c}
\sum\limits_{i \in B_{x/s}} {x^{(i)} \over s^{(i)}} \; \refGE{eq-XSLow1} \; (1-\beta) \rho(w) \sum\limits_{i \in B_{x/s}} {1 \over (s^{(i)})^2}\\
\refGE{eq-Gap} \; (1-\beta) \rho(w) \min\limits_{s \geq 0} \left\{ \sum\limits_{i \in B_{x/s}} {1 \over (s^{(i)})^2}: \; \sum\limits_{i \in B_{x/s}} s^{(i)} \leq {1 \over x^*_{\min}} \la s, x \ra \right\}
\\
= \; (1-\beta) \rho(w) {m^3 \over \la s, x \ra^2 } (x^*_{\min})^2 \; \geq \; (1-\beta) {m^3 (v^{(0)} - \| v \|^2)\over (n+1) (v^{(0)})^2 } (x^*_{\min})^2 \; \refGE{eq-PFin} \; m^3 {x^*_{\min} \over s^*_{\min}}.
\ea
$$
It remains to note that
$$
\ba{rcl}
\sum\limits_{i \not\in B_{x/s}} {x^{(i)} \over s^{(i)}} & \refLE{eq-XSLow1} & {1 \over (1-\beta) \rho(w)} \sum\limits_{i \not\in B_{x/s}} (x^{(i)})^2 \; \leq \; {1 \over (1-\beta) \rho(w)} \left( \sum\limits_{i \not\in B_{x/s}} x^{(i)} \right)^2\\
\\
& \refLE{eq-Gap} &  {\la s,x \ra^2 \over (1-\beta) \rho(w) (s^*_{\min})^2} \; \leq \; {(n+1) \mu^*(w) \over (1-\beta) (s^*_{\min})^2} \; \refLE{eq-PFin} \; {x^*_{\min} \over s^*_{\min}}.
\ea
$$
Thus, we get (\ref{eq-AwakeXS}).
\qed

The numerical verification of inequalities (\ref{eq-AwakeX}) - (\ref{eq-AwakeXS}) is very cheap. Therefore, in practical implementations of the parabolic target-following schemes, they can serve as conditions for activating the optimality test (\ref{eq-Cand}).

Straightforward implementation of the test (\ref{eq-Cand}) needs inversion of a non-symmetric matrix $A_{B_a}$. For the indicator ${x \over s}$, the cost of this operation can be reduced. Indeed, for the basis $B = B_{x/s}$, let us form the matrix $\Sigma_{x/s} = A_B X_B S_B^{-1} A_B^T$. Note that this matrix is a part of the full matrix $\Sigma = A^T XS^{-1}A$, which is required for computing affine-scaling directions. Hence, computation of $\Sigma_{x/s}$ does not entail any additional cost. However, since $\Sigma_{x/s}^{-1} = A_B^{-T} S_B X_B^{-1} A_B^{-1}$, we can use this matrix for computing the candidate optimal solution (\ref{eq-Cand}):
\beq\label{eq-Cand1}
\ba{rcl}
x^*_B & = & X_B S_B^{-1} A_B^T \Sigma_{x/s}^{-1} b, \quad y^*_{x/s} \; = \; \Sigma_{x/s}^{-1} A_B X_B S_B^{-1} c_B.
\ea
\eeq
In this case, the main term in the cost of the optimality test corresponds to computing a Cholesky factorization of a symmetric $m \times m$-matrix (this is ${m^3 \over 6}$).

\section{Numerical experiments}\label{sc-Num}
\SetEQ

For our computational experiments, we use a simple random generator proposed in \cite{AFF}. It works as follows.
\begin{itemize}
\item
Firstly, we generate a strictly feasible primal-dual pair of points $(\hat x, \hat s)$ for validating condition~(\ref{eq-SFeas}). Their entries are uniformly distributed in the interval $(0,1)$.
\item
After that, we form matrix $A \in \R^{m\times n}$ with entries  uniformly distributed in $(-1,1)$.
\item
Now, we can define $b = A \hat x$ and $c = \hat s$.
\item
The starting point $u_0$ for our methods is chosen as $(\hat x, \hat s, 0)$.
\end{itemize}

In the table below, we present preliminary computational results for random problems of small and medium dimensions with $32 \leq m \leq {n \over 2}$ and $64 \leq n \leq 1024$.
\beq\label{tab-Res}
\ba{|c|c|c|c|c|c|}
\hline
M \backslash N & 64 & 128 & 256 & 512 & 1024 \\
\hline
32 & 9.1 \pm 13.7 \% & 10.5 \pm 10.6 \% & 11.2 \pm 10.8 & 12.1 \pm 9.9\%  & 13.0 \pm 10.3 \% \\ \hline 
64 & & 11.4 \pm 11.1\% & 12.6 \pm 9.7\% &  13.7 \pm 8.6\% & 14.4 \pm 7.7\% \\ \hline
128 &  & & 13.6 \pm  8.6\% & 15.2 \pm 9.0\% & 16.2 \pm 7.5\% \\ \hline
256 & & & & 16.3 \pm 7.1\% & 18.0 \pm 7.3\% \\ \hline
512 & & & & & 19.2 \pm 7.4 \% \\
\hline
\ea
\eeq
In each cell, we put the average number of predictor steps of method (\ref{met-PT2}) required for reaching the accuracy $\epsilon = 10^{-8}$ in the duality gap. Our results correspond to the series of random test problems of length one hundred. The second value in the cell is the relative standard deviation in the series. We do not display the results for method (\ref{met-AutoPT}) since they are very similar to the results of method (\ref{met-AffinePT}), presented in \cite{AFF}. However, the performance of the second-order scheme (\ref{met-PT2}) appears to be much better. For the latter scheme, the required number of iterations is usually in 1.5 times smaller than that of methods (\ref{met-AffinePT}) or (\ref{met-AutoPT}).

In our opinion, these results are very promising. As in  numerical testing of \cite{AFF}, in {\em all} our experiments, each predictor step is followed by a {\em single} corrector step (hence, we do not display their counting). A quite accurate estimate of the number for predictor steps in method (\ref{met-PT2}) is given by the model
\beq\label{eq-For}
\ba{rcl}
k & \approx  & {1 \over 4}\left(25 + \log_2 m \cdot \log_2 {n \over 16} \right).
\ea
\eeq
In our experiments, the standard deviation of this forecast is $0.46$ iterations. We do not specify in this expression a  dependence on accuracy $\epsilon>0$ since for all our test problems method (\ref{met-PT2}) demonstrates an extremely fast local convergence. Typically, it goes even beyond the quadratic rate, as it was predicted by (\ref{eq-Rate3}).

The above numerical results serve as a serious motivation for testing the possible advantages of finite termination technique (see Section \ref{sc-Fin}) as applied to method (\ref{met-PT2}).
We present below our computational results for the indicator $xs^{-1}$. In the Table (\ref{tab-FRes}), the index for the average number of iterations shows how many problems in the whole series of 100 problems were terminated by the Termination Test (\ref{eq-Cand}). We accept there a real number $r$ to be non-negative if $r \geq - {\epsilon \over 100}$.
\beq\label{tab-FRes}
\ba{|c|c|c|c|c|c|}
\hline
M \backslash N & 64 & 128 & 256 & 512 & 1024 \\
\hline
32 & 7.0_{100} \pm 25 \% & 8.1_{99} \pm 16 \% & 8.6_{100} \pm 18\% & 9.6_{99} \pm 15\%  & 10.5_{98} \pm 15 \% \\ \hline 
64 & & 9.6_{85} \pm 22\% & 10.3_{87} \pm 16\% &  11.5_{98} \pm 14\% & 12.2_{94} \pm 14\% \\ \hline
128 &  & & 12.3_{64} \pm  16\% & 13.8_{67} \pm 14\% & 14.8_{77} \pm 13\% \\ \hline
256 & & & & 16.1_{15} \pm 8\% & 17.9_{13} \pm 8\% \\ \hline
512 & & & & & 19.2_0 \pm 7 \% \\
\hline
\ea
\eeq

As we can see, for small problems, the Termination Test works very well. However, when the dimensions increase, the fast local convergence becomes more and more important. For the biggest dimensions, the method almost always stops before the optimal basis could be detected by our tests. 

In all our experiments, the indicator
$$
\ba{rcl}
\beta(x,s) & = & {1 \over m^3} \sum\limits_{i \in B_{x/s}} {x^{(i)} \over s^{(i)}} \Big[ 
\sum\limits_{i \not\in B_{x/s}} {x^{(i)} \over s^{(i)}}\Big]^{-1}
\ea
$$
becomes big only in a couple of iterations before termination of the process (see (\ref{eq-AwakeXS})). Hence, the inequality $\beta(x,s) \geq 1$ can be used as an efficient activating condition for an attempt to guess the optimal primal basis.

\end{document}